\documentclass[11pt,a4paper,final, dvipsnames]{amsart}

\input xy
\xyoption{all}

\DeclareMathAlphabet{\mathpzc}{OT1}{pzc}{m}{it}

\usepackage{amsfonts,amsmath,amssymb,indentfirst,mathrsfs,amscd}
\usepackage[mathscr]{eucal}
\usepackage{graphicx}
\usepackage{nicefrac}

\usepackage{color}
\usepackage[utf8]{inputenc}
\usepackage[spanish, english]{babel}
\usepackage{booktabs}
\usepackage{multirow}
\usepackage{verbatim}

\usepackage{calrsfs}
\DeclareMathAlphabet{\omscal}{OMS}{zplm}{m}{n}

\usepackage[all]{xy}
\usepackage{graphicx}
\usepackage{stmaryrd}
\usepackage{enumitem}

\usepackage{empheq}
\usepackage{color}
\usepackage{xcolor}

\author[]{\'Angel Gonz\'alez-Prieto, Marina Logares and Vicente Mu\~noz}
\address{Escuela T\'ecnica Superior de Ingenieros de Sistemas Inform\'aticos, Universidad Polit\'ecnica de Madrid, Calle Alan Turing s/n (Carretera de Valencia Km 7), 28031 Madrid, Spain}\email{angel.gonzalez.prieto@upm.es}
\address{School of Computing, Electronics and Mathematics (Faculty of Science and Engineering), Plymouth University, 2-5, Kirkby Place,
Plymouth, PL1 8AA United Kingdom.}\email{marina.logares@plymouth.ac.uk}
\address{Departamento de \'Algebra, Geometr\'ia y Topolog\'ia, Facultad de Ciencias, Universidad de M\'alaga, Campus de Teatinos s/n, 29071 Málaga, Spain}\email{vicente.munoz@uma.es}

\title[]{A Lax Monoidal Topological Quantum Field Theory\\for Representation Varieties}


\DeclareMathOperator{\Hom}{Hom\,}           %Hom%
\DeclareMathOperator{\Gr}{Gr}

\usepackage{amsmath}
\begin{document}

\newtheorem{thm}{Theorem}[section]
\newtheorem{prop}[thm]{Proposition}
\newtheorem{lem}[thm]{Lemma}
\newtheorem{cor}[thm]{Corollary}
\newtheorem{conjecture}{Conjecture}

\theoremstyle{definition}
\newtheorem{defn}[thm]{Definition}
\newtheorem{ex}[thm]{Example}
\newtheorem{as}{Assumption}

\theoremstyle{remark}
\newtheorem{rmk}[thm]{Remark}

\theoremstyle{remark}
\newtheorem*{prf}{Proof}

\newcommand{\iacute}{\'{\i}} %i con acento%
\newcommand{\norm}[1]{\lVert#1\rVert} %norma%

\newcommand{\lto}{\longrightarrow}
\newcommand{\hra}{\hookrightarrow}

\newcommand{\suchthat}{\;\;|\;\;}
\newcommand{\dbar}{\overline{\partial}}

\newcommand{\cA}{\omscal{A}}
\newcommand{\cC}{\omscal{C}}
\newcommand{\cD}{\omscal{D}} %D-modules
\newcommand{\cF}{\omscal{F}} %Sheaf
\newcommand{\cG}{\omscal{G}} %Gauge group%
\newcommand{\cO}{\omscal{O}} %Holomorphic functions sheaf%
\newcommand{\cL}{\omscal{L}} %Lie derivative 
\newcommand{\cM}{\omscal{M}} %Moduli space% %moduli of parabolic bundles% %moduli of U(p,q) bundles%
\newcommand{\cN}{\omscal{N}} %Space of minimal points of the Morse function%
\newcommand{\cP}{\omscal{P}} %Moduli of K(D) pairs%
\newcommand{\cS}{\omscal{S}} %Moduli of solutions of Hitchin's equations, contructed by Konno%
\newcommand{\cU}{\omscal{U}} %Moduli of stable U(p,q) parabolic Higgs bundles%
\newcommand{\cX}{\omscal{X}}
\newcommand{\cZ}{\omscal{Z}}
\newcommand{\cT}{\omscal{T}}
\newcommand{\cV}{\omscal{V}}
\newcommand{\cB}{\omscal{B}}
\newcommand{\cR}{\omscal{R}}
\newcommand{\cH}{\omscal{H}}

\newcommand{\ext}{\mathrm{ext}} % an extension%
\newcommand{\x}{\times}

\newcommand{\mM}{\mathscr{M}} %Meromorphic function sheaf%

\newcommand{\CC}{\mathbb{C}} %Complex numbers%
\newcommand{\QQ}{\mathbb{Q}} %Rational numbers%
\newcommand{\PP}{\mathbb{P}} %projective space%
\newcommand{\HH}{\mathbb{H}} %Hypercohomology, quaternions..%
\newcommand{\RR}{\mathbb{R}} %Real numbers%
\newcommand{\ZZ}{\mathbb{Z}} %Integer numbers%
\newcommand{\NN}{\mathbb{N}} %Natural numbers%
\newcommand{\TT}{\mathbb{T}} %Torus numbers%

\newcommand{\acts}{\circlearrowright} %Torus numbers%
\newcommand{\Stab}{\mathrm{Stab}} %Stabilizer%

\renewcommand{\lg}{\mathfrak{g}} %Lie algebra of G%
\newcommand{\lh}{\mathfrak{h}} %Lie algebra of H%
\newcommand{\lu}{\mathfrak{u}} %Lie algebra of U%
\newcommand{\la}{\mathfrak{a}} %Lie algebra of A%
\newcommand{\lb}{\mathfrak{b}} %Lie algebra of B%
\newcommand{\lm}{\mathfrak{m}} %Lie algebra of M%
\newcommand{\lgl}{\mathfrak{gl}} %Lie algebra of GL%
\newcommand{\too}{\longrightarrow}
\newcommand{\imat}{\sqrt{-1}} %i%

% NEW COMMANDS

\newcommand\id{\mathrm{id}}

% Categories
\newcommand\Set{\mathbf{Set}}
\newcommand\CBord[3]{\mathbf{Bd}_{{#1 #3}}^{#2}}
\newcommand\CBordp[1]{\CBord{#1}{}{}}
\newcommand\CBordpp[1]{\mathbf{Bdp}_{{#1}}}
\newcommand\ClBordp[1]{\mathbf{l}\CBord{#1}{}{}}
\newcommand\CTub[3]{\mathbf{Tb}_{{#1 #3}}^{#2}}
\newcommand\CTubp[1]{\CTub{#1}{}{}}
\newcommand\CTubpp[1]{\mathbf{Tbp}_{#1}{}{}}
\newcommand\CTubppo[1]{\mathbf{Tbp}_{#1}^0}
\newcommand\CClose[3]{\mathbf{Cl}_{{#1 #3}}^{#2}}
\newcommand\CClosep[1]{\CClose{#1}{}{}}
\newcommand\Obj[1]{\mathrm{Obj}(#1)}
\newcommand\Mor[1]{\mathrm{Mor}(#1)}
\newcommand\Vect[1]{{#1}\textrm{-}\mathbf{Vect}}
\newcommand\Mod[1]{{#1}\textrm{-}\mathbf{Mod}}
\newcommand\Rng{\mathbf{Ring}}
\newcommand\Grpd{\mathbf{Grpd}}
\newcommand\Grpdo{\mathbf{Grpd}_0}
\newcommand\MHS[1]{\mathbf{MHS}}
\newcommand\CVar{\mathbf{Var}_\CC}
\newcommand\PHM[2]{\cM_{#1}^p(#2)}
\newcommand\MHM[1]{\cM_{#1}}
\newcommand\VMHS[1]{VMHS({#1})}
\newcommand\geoVMHS[1]{VMHS_g({#1})}
\newcommand\Perv[2]{\mathrm{Perv}({#1},{#2})}
\newcommand\Par[1]{\mathrm{Par}({#1})}
\newcommand\Dcs[2]{D^b_{cs}({#1},{#2})}
\newcommand\K[1]{\mathrm{K}#1}
\newcommand\KQ[1]{\omscal{K} #1}
\newcommand\Ab{\mathbf{Ab}}
\newcommand\VarOver[1]{\mathrm{Var}/{#1}}
\newcommand\VarOverH[1]{\mathrm{Var}_H/{#1}}
\newcommand\CPP{\cP\cP}
\newcommand\Bim[1]{{#1}\textrm{-}\mathbf{Bim}}
\newcommand\Span[1]{\mathrm{Span}({#1})}
\newcommand\Spano[1]{\mathrm{Span}^{op}({#1})}
\newcommand\Geo{\cG}
\newcommand\AlgA{\cA}
\newcommand\AlgB{\cB}

% Hodge theory
\newcommand\RDelHod{e}
\newcommand\DelHod[1]{e\left(#1\right)}
\newcommand\eVect{\omscal{E}}
\newcommand\e[1]{\eVect\left(#1\right)}
\newcommand\intMor[2]{\int_{#1}\,#2}

\newcommand{\rat}{\mathrm{rat}}

% Representation Varieties
\newcommand\Dom[1]{\omscal{D}_{#1}}

\newcommand\Xf[1]{\omscal{X}_{#1}}					% Free non-parabolic
\newcommand\Xs[1]{X_{#1}}							% Surface group non-parabolic

\newcommand\Xft[2]{\overline{\omscal{X}}_{#1, #2}} % Free tr=2
\newcommand\Xst[2]{\overline{X}_{#1, #2}}			% Surface group tr=2

\newcommand\Xfp[2]{\omscal{X}_{#1, #2}}			% Free J+
\newcommand\Xsp[2]{X_{#1, #2}}						% Surface group J+

\newcommand\Xfd[2]{\omscal{X}_{#1; #2}}			% Free diagonal
\newcommand\Xsd[2]{X_{#1; #2}}						% Surface group diagonal

\newcommand\Xfm[3]{\omscal{X}_{#1, #2; #3}}		% Free mixed
\newcommand\Xsm[3]{X_{#1, #2; #3}}					% Surface group mixed

\newcommand\XD[1]{#1^D}
\newcommand\XDh[1]{#1^{\to D}}
\newcommand\XU[1]{#1^U}
\newcommand\XUh[1]{#1^{\to U}}
\newcommand\XP[1]{#1^{\pm 1}}
\newcommand\XPh[1]{#1^{\to \pm 1}}
\newcommand\XI[1]{#1^{I}}
\newcommand\XTilde[1]{#1^{R}}
\newcommand\Xred[1]{#1^{red}}
\newcommand\Xirred[1]{#1^{irred}}
\newcommand\Rep[1]{\mathfrak{X}_{#1}}

% Miscelany
\newcommand\RM[2]{R\left(\left.#1\right|#2\right)}
\newcommand\RMc[3]{R_{#1}\left(\left.#2\right|#3\right)}
\newcommand\set[1]{\left\{#1\right\}}

\newcommand\sZ{\mathcal{Z}}
\newcommand\lZ{\cZ}
\newcommand\Aff[1]{\textrm{AGL}(#1)}
\newcommand\ASO[1]{\textrm{ASO}(#1)}
\newcommand\GL{\textrm{GL}}
\newcommand\PGL{\textrm{PGL}}
\newcommand\SL{\textrm{SL}}

\hyphenation{mul-ti-pli-ci-ty}

\hyphenation{mo-du-li}

\begin{abstract}
We construct a lax monoidal Topological Quantum Field Theory that computes Deligne-Hodge polynomials of representation varieties of the fundamental group of any closed manifold into any complex algebraic group $G$. As byproduct, we obtain formulas for these polynomials in terms of homomorphisms between the space of mixed Hodge modules on $G$. The construction is developed in a categorical-theoretic framework allowing its application to other situations.
\end{abstract}
\null
\vspace{-1.1cm}
\maketitle

\setcounter{section}{0}

\vspace{-1cm}

%%%%%%%%%%%%%%%%%%%%%%%%%%%%%%%%%%%%%%%%%%%%%%%%%%%%%%%%%%%%%%%%
%%%%%%%%%%%%%%%%%%% SECTION: INTRODUCTION %%%%%%%%%%%%%%%%%%%%%%
%%%%%%%%%%%%%%%%%%%%%%%%%%%%%%%%%%%%%%%%%%%%%%%%%%%%%%%%%%%%%%%%

\section{Introduction}
\let\thefootnote\relax\footnotetext{\noindent \emph{2010 Mathematics Subject Classification}. Primary:
 57R56. % TQFT
 Secondary:
 14C30, %Hodge theory
 14D07, % Variations of Hodge structures
 14D21. % Applications of moduli in mathematical physics
 
\emph{Key words and phrases}: TQFT, moduli spaces, $E$-polynomial, representation varieties.}

Let $W$ be a compact manifold, possibly with boundary, and let $G$ be a complex algebraic group. The set of representations $\rho: \pi_1(W) \to G$ can be endowed with a complex algebraic variety structure, the so-called representation variety of $W$ into $G$, denoted $\Rep{G}(W)$. The group $G$ itself acts on $\Rep{G}(W)$ by conjugation so, taking the Geometric Invariant Theory (GIT) quotient of $\Rep{G}(W)$ by this action (see \cite{Newstead:1978}) we obtain $\cM_G(W) = \Rep{G}(W) \sslash G$, the moduli space of representations of $\pi_1(W)$ into $G$, as treated in \cite{Nakamoto}. It is customary to call these spaces character varieties or, in the context of non-abelian Hodge theory, Betti moduli spaces. Even in the simplest cases $G=\GL(n, \CC), SL(n,\CC)$ and $W=\Sigma$, a closed orientable surface, the topology of these varieties is extremely rich and has been the object of studies during the past twenty years.

One of the main reasons of this study is the prominent role of these varieties in the non-abelian Hodge theory. For $G = \GL(n, \CC)$ (resp.\ $G = \SL(n, \CC)$), an element of $\cM_G(\Sigma)$ defines a $G$-local system and, thus, a rank $n$ algebraic bundle $E \to \Sigma$ of degree $0$ (resp.\ and fixed determinant) with a flat connection $\nabla$ on it. Hence, the Riemann-Hilbert correspondence (\cite{SimpsonI} \cite{SimpsonII}) gives a real analytic correspondence between $\cM_G(\Sigma)$ and the moduli space of flat bundles of rank $n$ and degree $0$ (and fixed determinant if $G=\SL(n,\CC)$), usually called the de Rham moduli space.

Furthermore, via the Hitchin-Kobayashi correspondence (\cite{Simpson:1992} \cite{Corlette:1988}), we also have that, for $\Sigma$ a compact Riemann surface and $G=\GL(n, \CC)$ (resp.\ $G=\SL(n,\CC)$), the Betti moduli space $\cM_G(\Sigma)$ is real analytic equivalent to the Dolbeault moduli space, that is, the moduli space of rank $n$ and degree $0$ (resp.\ and fixed determinant) $G$-Higgs bundles i.e.\ bundles $E \to \Sigma$ together with a field $\Phi: E \to E \otimes K_\Sigma$ called the Higgs field.

Motivated by these correspondences, we can also consider representation varieties of a manifold $W$ with a parabolic structure $Q$. This $Q$ consists of a finite set of pairwise disjoint subvarieties $S_1, \ldots, S_r$ of $W$ of codimension $2$ and conjugacy classes $\lambda_1, \ldots, \lambda_r \subseteq G$ that allows us to define $\Rep{G}(W, Q)$ as the set of representations $\rho: \pi_1(W - S_1 - \ldots - S_r) \to G$ such that the image of the loops arround $S_i$ must live in $\lambda_i$ (see section \ref{subsec:parabolic-case} for precise definition). Analogously, we set $\cM_G(W, Q)= \Rep{G}(W,Q) \sslash G$.

When $W=\Sigma$ is a surface, the $S_i$ are a set of (marked) points, called the parabolic points, and we can obtain stronger results. For example, for $G=\GL(n,\CC)$ and $Q$ a single marked point $p \in \Sigma$ and $\lambda = \left\{e^{2\pi i d/n}\,Id\right\}$ we have that $\cM_G(\Sigma, Q)$ is diffeomorphic to the moduli space of rank $n$ and degree $d$ Higgs bundles and to the moduli space of rank $n$ logarithmic flat bundles of degree $n$ with a pole at $p$ with residue $-\frac{d}{n}Id$. In this case, $\cM_G(\Sigma, Q)$ is referred to as the twisted caracter variety.

For an arbitrary number of marked points $p_1, \ldots, p_r \in \Sigma$ and different semi-simple conjugacy classes of $G=\GL(n,\CC)$, we obtain diffeomorphisms between moduli spaces of parabolic Higgs bundles with parabolic structures (with general weights) on $p_1, \ldots, p_r$ and with the moduli space of logarithmic flat connections with poles on $p_1, \ldots, p_r$ (\cite{Simpson:parabolic}). Incidentally, other correspondences can also appear as for the case of $G=\SL(2,\CC)$, $\Sigma$ an elliptic curve and $Q$ two marked points with different semi-simple conjugacy classes not containing a multiple of identity, in which $\cM_G(\Sigma, Q)$ is diffeomorphic to the moduli space of doubly periodic instantons through the Nahm transform \cite{Biquard-Jardim} \cite{Jardim}.

Using these correspondences, it is possible to compute the Poincar\'e polynomial of character varieties by means of Morse theory. Following these ideas, Hitchin, in the seminal paper \cite{Hitchin}, gave the Poincar\'e polynomial for $G=\SL(2,\CC)$ in the non-parabolic case, Gothen also computed it for $G=\SL(3,\CC)$ in \cite{Gothen} and Garc\'ia-Prada, Heinloth and Schmitt for $G=\GL(4,\CC)$ in \cite{GP-Heinloth-Schmitt}. In general, in \cite{Schiffmann:2016} and \cite{Mozgovoy-Schiffmann} (see also \cite{Mellit}) a combinatorial formula is given for arbitrary $G=\GL(r, \CC)$ provided that $n$ and $d$ are coprime. In the parabolic case, Boden and Yokogawa calculated it in \cite{Boden-Yokogawa} for $G=\SL(2,\CC)$ and generic semi-simple conjugacy classes and Garc\'ia-Prada, Gothen and the third author for $G=\GL(3,\CC)$ and $G=\SL(3,\CC)$ in \cite{GP-Gothen-Munoz}. 

However, these correspondences from non-abelian Hodge theory are far from being algebraic. Hence, the study of their (mixed) Hodge structure on cohomology turns important. A useful combinatorial tool for this purpose is the so-called Deligne-Hodge polynomial, also referred to as $E$-polynomial, that, to any complex algebraic $X$ assigns a polynomial $\DelHod{X} \in \ZZ[u,v]$. As described in section \ref{subsec:mhs}, this polynomial is constructed as an alternating sum of the Hodge numbers of $X$, in the spirit of a combination of Poincar\'e polynomial and Euler characteristic.

A great effort has been made to compute these $E$-polynomials for character varieties. The first strategy was accomplished by \cite{Hausel-Rodriguez-Villegas:2008} by means of a theorem of Katz of arithmetic flavour based on the Weil conjectures and the Lefschetz principle. Following this method, when $\Sigma$ is an orientable surface, an expression of the $E$-polynomial for the twisted character varieties is given in terms of generating functions in \cite{Hausel-Rodriguez-Villegas:2008} for $G=\GL(n,\CC)$ and in \cite{Mereb} for $G=\SL(n,\CC)$. Recently, using this technique, explicit expressions of the $E$-polynomials have been computed, in \cite{Baraglia-Hekmati:2016}, for the untwisted case and orientable surfaces with $G=\GL(3,\CC),\SL(3,\CC)$ and for non-orientable surfaces with $G=\GL(2,\CC),\SL(2,\CC)$. Also they have checked the formulas given in \cite{MM} for orientable surfaces and $G=\SL(2,\CC)$ in the untwisted case. 

The other approach to this problem was initiated by the second and third authors together with Newstead in \cite{LMN}. In this case, the strategy is to focus on the computation of $\DelHod{\Rep{G}(\Sigma)}$ and, once done, to pass to the quotient. In this method, the representation variety is chopped into simpler strata for which the $E$-polynomial can be computed. After that, one uses the additivity of $E$-polynomials to combine them and get the one of the whole space.

Following this idea, in the case $G=\SL(2,\CC)$, they computed the $E$-polynomials of character varieties for a single marked point and genus $g=1,2$ in \cite{LMN}. Later, the second and third authors computed them for two marked points and $g=1$ in \cite{LM} and the third author with Mart\'inez for a marked point and $g=3$ in \cite{MM:2016}. In the case of arbitrary genus and, at most, a marked point, the case $G=\SL(2,\CC)$ was accomplished in \cite{MM} and the case $G=\PGL(2,\CC)$ in \cite{Martinez:2017}.

In these later papers, this method is used to obtain recursive formulas of $E$-polynomials of representation varieties in terms of the ones for smaller genus. This recursive nature for character varieties is widely present in the literature, as in and \cite{Mozgovoy:2012}, \cite{Hausel-Letellier-Villegas:2013}, \cite{Diaconescu:2017} and \cite{Carlsson-Rodriguez-Villegas}. It suggests that some sort of recursion formalism, in the spirit of Topological Quantum Field Theory (TQFT for short), must hold. That is the starting point of the present paper.

In the parabolic case, much remains to be known. The most important advance was given in \cite{Hausel-Letellier-Villegas}, following the arithmetic method, for $G=\GL(n,\CC)$ and generic semi-simple marked points. Using the geometric method, as we mentioned above, only at most two marked points have been studied.

In this paper, we propose a general framework for studing $E$-polynomials of representation varieties based on the stratification strategy, valid for any complex algebraic group $G$, any manifold (not necessarely surfaces) and any parabolic configuration. For this purpose, section \ref{sec:hodge-theory} is devoted to review the fundamentals of Hodge theory and Saito's mixed Hodge modules as a way of tracing variations of Hodge structures with nice functorial properties (see \cite{Saito:1990}).

With these tools at hand, we can develop a categorical theoretic machinery that shows how recursive computations of Deligne-Hodge polynomials can be accomplished. Based on TQFTs (i.e.\ monoidal functors $Z: \CBordp{n} \to \Vect{k}$ between the category of $n$-bordisms and the category of $k$-vector spaces, as introduced in \cite{Atiyah:1988}) we define, in section \ref{subsec:bimodules}, a weaker version of them in the context of $2$-categories and pairs of spaces. We propose to consider lax monoidal lax functors $\sZ: \CBordpp{n} \to \Bim{R}$ between the $2$-category of pairs of bordisms and the $2$-category of $R$-algebras and bimodules (being $R$ a ring), which we called soft Topological Quantum Field Theories of pairs. These soft TQFTs are, in some sense, parallel to the so-called Extended Topological Quantum Field Theories, as studied in \cite{Baez-Dolan}, \cite{Freed:1994} or \cite{Lauda-Pfeiffer} amongst others.

In this setting, we will show how a soft TQFT can be constructed in full generality from two basic pieces: a functor $\Geo: \CBordpp{n} \to \Span{\CVar}$ (being $\Span{\CVar}$ the $2$-category of spans of the category of complex algebraic varieties), called the geometrisation, and a contravariant functor $\AlgA: \CVar \to \Rng$, called the algebraisation. In section \ref{sec:LTQFT-E-pol}, we will apply these ideas to the computation of Deligne-Hodge polynomials of representations varieties. For that, we will define the geometrisation by means of the fundamental groupoid of the underlying manifold and we will use the previously developed theory of mixed Hodge modules for an algebraisation. 

Even though this soft TQFT encodes the recursive nature of the Deligne-Hodge polynomial, we can make it even more explicit. For this purpose, in section \ref{subsec:modules-twists}, we define a lax monoidal Topological Quantum Field Theories of pairs as a lax monoidal strict functor $\lZ: \CBordpp{n} \to \Mod{R}_t$, where $\Mod{R}_t$ is the usual category of $R$-modules with an additional $2$-category structure (see definition \ref{defn:twisted-modules}). In this context, we show how a partner covariant functor $\AlgB: \CVar \to \Mod{R}$ to the algebraisation $\AlgA$ allow us to define a natural lax monoidal TQFT. Again, we will use this idea to define a TQFT computing $E$-polynomials of representation varieties. In this formulation, an explicit formula for these polynomials can be deduced.

\begin{thm}\label{thm:existence-s-LTQFT}
There exists a lax monoidal TQFT of pairs, $\lZ: \CBordpp{n} \to \Mod{R}_t$, where $R=\K{\MHS{\QQ}}$ is the $K$-theory ring of the category of mixed Hodge modules, such that, for any $n$-dimensional connected closed orientable manifold $W$ and any non-empty finite subset $A \subseteq W$ we have
$$
	\DelHod{\lZ(W, A)(\QQ_0)} = \DelHod{G}^{|A|-1} \DelHod{\Rep{G}(W)}.
$$ 
\end{thm}

It is well worthy to point out that the construction of this lax monoidal TQFT follows the underlying philosophy of quantifying the desired invariant for bordisms and computing the associated homomorphism as a pullback from the ingoing boundary to the bordism followed by a pushout onto the outgoing boundary, following the so-called 'push-pull construction' (see, for example, \cite{Freed:1994}, \cite{Freed-Lurie:2010}, \cite{Freed-Hopkins-Teleman:2010}, \cite{Haugseng}, \cite{Ben-Zvi-Nadler:2016} or \cite{Ben-Zvi-Gunningham-Nadler} among others).

This idea of using a TQFT for understanding the cohomology of character varieties has been successfully applied several times, as in \cite{Ben-Zvi-Nadler}, \cite{Ben-Zvi-Francis-Nadler:2010}, \cite{Kassabov-Patotski} and in the recent paper \cite{Ben-Zvi-Gunningham-Nadler}. The exact interplay between these constructions must be addressed in future work. In particular, it could be expectable a strong relation between the construction of \cite{Ben-Zvi-Gunningham-Nadler}, based on the stack of G-local systems and Lurie’s cobordisms hypothesis, and the one given in this paper. 

With a view towards applications, in section \ref{subsec:almost-TQFT}, we will show how, for computational purposes, it is enough to consider tubes instead of general bordisms, defining what we call an almost-TQFT, $\mathfrak{Z}: \CTubpp{n} \to \Mod{R}$. The idea is that, given a closed surface $\Sigma$, we can choose a suitable handlebody decomposition of $\Sigma$ as composition of tubes. Therefore, in order to compute $\lZ(\Sigma)$ we do not need the knowledge of the image of any general cobordism but just of a few tubes.

For $n=2$, we will give explicitly the image of the generators of $\CTubpp{2}$ for the corresponding almost-TQFT computing $E$-polynomials of representation varieties. From this description, we obtain an explicit formula of $\DelHod{\Rep{G}(\Sigma, Q)}$ in terms of simpler pieces, see theorem \ref{thm:almost-tqft-parabolic}. From this formula, a general algorithm for computing these polynomials arises. Actually, the computations of \cite{LMN} and \cite{MM} are particular calculations of that program.

This paper is part of the PhD Thesis \cite{Gonzalez-Prieto:Thesis} of the first author under the supervision of the second and third autors. In that thesis, the computational method developed in this paper is used for computing the Deligne-Hodge polynomials of $\SL(2,\CC)$-parabolic character varieties with punctures of Jordan type. This program is also contained in the preprints \cite{GP-2018a} and \cite{GP-2018b}. The case of punctures of semi-simple type is addressed in the upcoming paper \cite{GP-2019}.

Another framework in which character varieties are central is the Geometric Langlands program (see \cite{Beilinson-Drinfeld}). In this setting, the Hitchin fibration satisfies the Strominger-Yau-Zaslow conditions of mirror symmetry for Calabi-Yau manifolds (see \cite{Strominger-Yau-Zaslow}) from which arise several questions about relations between $E$-polynomials of character varieties for Langland dual groups as conjectured in \cite{Hausel:2005} and \cite{Hausel-Rodriguez-Villegas:2008}. The validity of these conjectures has been discussed in some cases as in \cite{LMN} and \cite{Martinez:2017}. Despite of that, the general case remains unsolved.
We hope that the ideas introduced in this paper could be useful to shed light into these questions.

\subsection*{Acknowledgements}

The authors want to thank Bruce Bartlett, David Ben-Zvi, Christopher Douglas and Constantin Teleman for useful conversations. We want to express our highest gratitude to Thomas Wasserman for his invaluable help throughout the development of this paper.

We thank the hospitality of the Mathematical Institute at University of Oxford offered during a research visit which was supported by the Marie Sklodowska Curie grant GREAT - DLV-654490. The work of the first and third authors has been partially supported by MINECO (Spain) Project MTM2015-63612-P. The first author was also supported by a "\!la Caixa" scholarship for PhD studies in Spanish Universities from "\!la Caixa" Foundation. The second author was supported by the Marie Sklodowska Curie grant GREAT - DLV-654490.

%%%%%%%%%%%%%%%%%%%%%%%%%%%%%%%%%%%%%%%%%%%%%%%%%%%%%%%%%%%%%%%%
%%%%%%%%%%%%%%%%%% SECTION: HODGE THEORY %%%%%%%%%%%%%%%%%%%%%%%
%%%%%%%%%%%%%%%%%%%%%%%%%%%%%%%%%%%%%%%%%%%%%%%%%%%%%%%%%%%%%%%%

\section{Hodge theory}
\label{sec:hodge-theory}

\subsection{Mixed Hodge structures}
\label{subsec:mhs}

Let $X$ be a complex algebraic variety. The rational cohomology of $X$, $H^\bullet(X; \QQ)$, carries an additional linear structure, called \emph{mixed Hodge structure}, which generalize the so-called pure Hodge structures. For further information, see \cite{DeligneII:1971} and \cite{DeligneIII:1971}, also \cite{Peters-Steenbrink:2008}.

\begin{defn}
Let $V$ be a finite dimensional rational vector space and let $k \in \ZZ$. A \emph{pure Hodge structure} of weight $k$ on $V$ consists of a finite decreasing filtration $F^\bullet$ of $V_\CC = V \otimes_\QQ \CC$
$$
	V_\CC \supseteq \ldots \supseteq F^{p-1}\,V \supseteq F^p\,V \supseteq F^{p+1}\,V \supseteq \ldots \supseteq \left\{0\right\}
$$
such that $F^p\,V \oplus \, \overline{F^{k-p+1}\,V} = V_\CC$ where conjugation is taken with respect to the induced real structure.
\end{defn}

\begin{rmk}
\label{rmk:PHS}
An equivalent description of a pure Hodge structure is as a finite decomposition
$$
	V_\CC = \bigoplus_{p+q=k} V^{p,q}
$$
for some complex vector spaces $V^{p,q}$ such that $\overline{V^{p,q}} \cong V^{q,p}$ with respect to the natural real structure of $V_\CC$. From this description, the filtration $F^\bullet$ can be recovered by taking $F^p\,V = \bigoplus_{r \geq p} V^{r,k-r}$. In these terms, classical Hodge theory shows that the cohomology of a compact K\"ahler manifold $M$ carries a pure Hodge structure induced by Dolbeault cohomology by $H^k(M;\CC) = \bigoplus\limits_{p+q = k} H^{p,q}(M)$. See \cite{Peters-Steenbrink:2008} for further information.
\end{rmk}

\begin{ex}
Given $m \in \ZZ$, we define the \emph{Tate structures}, $\QQ(m)$, as the pure Hodge structure whose underlying rational vector space is $(2\pi i)^m \QQ \subseteq \CC$ with a single-piece decomposition $\QQ(m) = \QQ(m)^{-m,-m}$. Thus, $\QQ(m)$ is a pure Hodge structure of weight $-2m$. Moreover, if $V$ is another pure Hodge structure of weight $k$ then $V(m) := V \otimes \QQ(m)$ is a pure Hodge structure of weight $k-2m$, called the \emph{Tate twist} of $V$. For short, we will denote $\QQ_0 = \QQ(0)$, the Tate structure of weight $0$. Recall that there is a well defined tensor product of pure (and mixed) Hodge structures, see Examples 3.2 of \cite{Peters-Steenbrink:2008} for details.
\end{ex}

\begin{defn}
Let $V$ be a finite dimensional rational vector space. A \emph{(rational) mixed Hodge structure} on $V$ consist of a pair of filtrations:
\begin{itemize}
	\item An increasing finite filtration $W_\bullet$ of $V$, called the \emph{weight filtration}.
	\item A decreasing finite filtration $F^\bullet$ of $V_\CC$, called the \emph{Hodge filtration}.
\end{itemize}
Such that, for any $k \in \ZZ$, the induced filtration of $F^\bullet$ on the graded complex $(Gr^W_k\,V)_\CC = \left(\frac{W_k\,V}{W_{k-1}\,V}\right)_\CC$ gives a pure Hodge structure of weight $k$. Given two mixed Hodge structures $(V, F, W)$ and $(V', F', W')$, a morphism of mixed Hodge structures is a linear map $f: V \to V'$ preserving both filtrations.
\end{defn}

Deligne proved in \cite{DeligneII:1971} and \cite{DeligneIII:1971} (for a concise exposition see also \cite{Zein-Trang:2014}) that, if $X$ is a complex algebraic variety, then $H^k(X; \QQ)$ carries a mixed Hodge structure in a functorial way. More preciselly, let $\CVar$ be the category of complex varieties with morphisms given by the regular maps, $\Vect{\QQ}$ the category of $\QQ$-vector spaces and $\MHS{\QQ}$ be the category of mixed Hodge structures. First, we have that $\MHS{\QQ}$ is an abelian category (see Théorème 2.3.5 of \cite{DeligneII:1971}) and, moreover, the cohomology functor $H^k(-;\QQ): \CVar \to \Vect{\QQ}$ factorizes through $\MHS{\QQ}$, that is
\[
\begin{displaystyle}
   \xymatrix
   {
  	 \CVar \ar[r]^{\hspace*{-8pt}H^k(-; \QQ)} \ar[d] & \Vect{\QQ} \\
  	 \MHS{\QQ} \ar[ru]
   }
\end{displaystyle}   
\]

\begin{rmk}
A pure Hodge structure of weight $k$ is, in particular, a mixed Hodge structure by taking the weight filtration with a single step. When $X$ is a smooth complex projective variety, the induced pure Hodge structure given by Remark \ref{rmk:PHS} corresponds to the mixed Hodge structure given above.
\end{rmk}

An analogous statement holds for compactly supported cohomology, that is $H_c^k(X;\QQ)$ has a mixed Hodge structure in a functorial way (see section 5.5 of \cite{Peters-Steenbrink:2008} for a complete construction). From this algebraic structure, some new invariants can be defined (see Definition 3.1 of \cite{Peters-Steenbrink:2008}). Given a complex algebraic variety $X$, we define the $(p,q)$-pieces of its $k$-th compactly supported cohomology groups by
$$
	H^{k;p,q}_c(X) := \Gr_p^{F}\left(\Gr_W^{p+q}\,H^k_c(X;\QQ)\right)_\CC
$$
From them, we define the \emph{Hodge numbers} as $h^{k;p,q}_c(X) = \dim\,H^{k;p,q}_c(X)$ and the \emph{Deligne-Hodge polynomial}, or \emph{$E$-polynomial}, as the alternating sum
$$
	\DelHod{X} = \sum_k (-1)^k h_c^{k;p,q}(X)\;u^pv^q \in \ZZ[u^{\pm 1},v^{\pm 1}]
$$

\begin{rmk}
Sometimes in the literature, the $E$-polynomial is defined as $\DelHod{X}(-u,-v)$. It does not introduce any important difference but it would produce an annoying change of sign.
\end{rmk}

\begin{rmk}
An important fact is that the K\"unneth isomorphism
$$
	H_c^\bullet(X; \QQ) \otimes H_c^\bullet(Y; \QQ) \cong H_c^\bullet(X \times Y; \QQ) 
$$
is an isomorphism of mixed Hodge structures (see Theorem 5.44 in \cite{Peters-Steenbrink:2008}). In particular, this implies that $\DelHod{X \times Y} = \DelHod{X} \DelHod{Y}$. When, instead of product varieties, we consider general fibrations, the monodromy plays an important role (see, for example \cite{LMN}, \cite{MM:2016} or \cite{MM}). The best way to deal with this issue is through the theory of mixed Hodge modules (see next).
\end{rmk}

\subsection{Mixed Hodge modules}
\label{subsec:mixed-hodge-modules}

In \cite{Saito:1990} (see also \cite{Saito:1986} and \cite{Saito:1989}) Saito proved that we can assign, to every complex algebraic variety $X$, an abelian category $\MHM{X}$ called the category of mixed Hodge modules on $X$. As described in \cite{Schurmann:2011}, if $X$ is smooth, the basic elements of $\MHM{X}$ are tuples $M = (S, F^\bullet, W_\bullet, K, \alpha)$ where $S$ is a regular holonomic $D_X$-module with $F^\bullet$ a good filtration, $K$ is a perverse sheaf (sometimes also called a perverse complex) of rational vector spaces and $W_\bullet$ is a pair of increasing filtrations of $S$ and $K$. These filtrations have to correspond under the isomorphism
$$
	\alpha: DR_X(S) \stackrel{\cong}{\to} K \otimes_{\underline{\QQ}_X} \underline{\CC}_X
$$
where $\underline{\CC}_X, \underline{\QQ}_X$ are the respective constant sheaves on $X$ and $DR_X$ is the Riemann-Hilbert correspondence functor between the category of filtered $D_X$-modules and the category of rational perverse sheaves on $X$, $\Perv{X}{\QQ}$ (for all these concepts see \cite{Peters-Steenbrink:2008}).
Starting with these basic elements, the category $\MHM{X}$ is constructed as a sort of ``controlled" extension closure of these tuples, in the same spirit as mixed Hodge structures are a closure of pure Hodge structures under extension. In the case that $X$ is singular, the construction is similar but more involved using local embeddings of $X$ into manifolds (see \cite{Saito:1990}, \cite{Peters-Steenbrink:2008} or \cite{Schurmann:2011}).

In \cite{Saito:1986} and \cite{Saito:1990} Saito proves that $\MHM{X}$ is an abelian category endowed with a functor
$$
	\rat_X: \MHM{X} \to \Perv{X}{\QQ}
$$
that extends to the (bounded) derived category as a functor
$$
	\rat_X: D^b\MHM{X} \to \Dcs{X}{\QQ}
$$
where $\Dcs{X}{\QQ}$ is the derived category of cohomological constructible complexes of sheaves that contains $\Perv{X}{\QQ}$ as a full abelian subcategory (\cite{Peters-Steenbrink:2008}, Lemma 13.22). Moreover, given a regular morphism $f: X \to Y$, there are functors
$$
	f_*, f_!: D^b\MHM{X} \to D^b\MHM{Y} \hspace{1cm} f^*, f^!: D^b\MHM{Y} \to D^b\MHM{X}
$$
which lift to the analogous functors on the level of constructible sheaves. Finally, the tensor and external product of constructibles complexes lift to bifunctors
$$
	\otimes: D^b\MHM{X} \times D^b\MHM{X} \to D^b\MHM{X} \hspace{1cm} \boxtimes: D^b\MHM{X} \times D^b\MHM{Y} \to D^b\MHM{X \times Y}.
$$

\begin{rmk}
\begin{itemize}
	\item Recall that $f_*, f_!$ on $\Dcs{X}{\QQ}$ are just the usual direct image and proper direct image on sheaves, $f^*$ is the inverse image sheaf and $f^!$ is the adjoint functor of $f_!$, the so-called extraordinary pullback. See \cite{Peters-Steenbrink:2008}, Chapter 13, for a complete definition of these functors.
	\item As in the case of constructible complexes, the external product can be defined in terms of the usual tensor product by
	$$
		M^\bullet \boxtimes N^\bullet = \pi_1^* M^\bullet \otimes \pi_2^* N^\bullet
	$$
for $M^\bullet \in D^b\MHM{X}$, $N^\bullet \in D^b\MHM{Y}$ and $\pi_1: X \times Y \to X$, $\pi_2: X \times Y \to Y$ the corresponding projections. 
	\end{itemize}
\end{rmk}

A very important feature of these induced functors is that they behave in a functorial way, as the following result shows. The proof of this claim is a compendium of Proposition 4.3.2 and Section 4.4 (in particular 4.4.3) of \cite{Saito:1990}.

\begin{thm}[Saito]\label{thm:funct-induced-functors-mhm}
The induced functors commute with composition. More explicitly, let $f: X \to Y$ and $g: Y \to Z$ regular morphisms of complex algebraic varieties, then
$$
	(g \circ f)_* = g_* \circ f_* \hspace{1cm} (g \circ f)_! = g_! \circ f_! \hspace{1cm} (g \circ f)^* = f^* \circ g^* \hspace{1cm} (g \circ f)^! = f^! \circ g^!
$$
Furthermore, suppose that we have a cartesian square of complex algebraic varieties (i.e.\ a pullback diagram in $\CVar$)
\[
\begin{displaystyle}
   \xymatrix
   {
   		W \ar[r]^{g'} \ar[d]_{f'} & X \ar[d]^{f} \\
   		Y \ar[r]_{g} & Z
   }
\end{displaystyle}   
\]
Then we have a natural isomorphism of functors $g^* \circ f_! \cong f'_! \circ (g')^*$.
\end{thm}

Given a complex algebraic variety $X$, associated to the abelian category $\MHM{X}$, we can consider the Grothendieck group, also known as $K$-theory group, denoted by $\K{\MHM{X}}$. Recall that it is the free abelian group generated by the objects of $\MHM{X}$ quotiented by the relation $M \sim M' + M''$ if $0 \to M' \to M \to M'' \to 0$ is a short exact sequence in $\MHM{X}$. By definition, we have an arrow on objects $\MHM{X} \to \K{\MHM{X}}$. Moreover, given $M^\bullet \in D^b\MHM{X}$, we can associate to it the element of $\K{\MHM{X}}$
$$
	[M^\bullet] = \sum_k (-1)^k H^k(M^\bullet) \in \K{\MHM{X}}
$$
where $H^k(M^\bullet) \in \MHM{X}$ is the $k$-th cohomology of the complex. This gives an arrow on objects $D^b\MHM{X} \to \K{\MHM{X}}$. Under this arrow, tensor product $\otimes: D^b\MHM{X} \times D^b\MHM{X} \to D^b\MHM{X}$ descends to a bilinear map $\otimes: \K{\MHM{X}} \times \K{\MHM{X}} \to \K{\MHM{X}}$ that endows $\K{\MHM{X}}$ with a natural ring structure.

With respect to induced morphisms, given $f:X \to Y$, the functors $f_*, f_!, f^*, f^!$ of mixed Hodge modules also descend to give group homomorphisms $f_*, f_!: \K{\MHM{X}} \to \K{\MHM{Y}}$ and $f^*, f^!: \K{\MHM{Y}} \to \K{\MHM{X}}$ (see Section 4.2 of \cite{Schurmann:2011}). For example, we define $f_!: \K{\MHM{X}} \to \K{\MHM{Y}}$ by
$$
	f_![M] := [f_!M] = \sum_k (-1)^k \,H^k\left(f_!M\right)
$$
where $[M]$ denotes the class of $M \in \MHM{X}$ on $\K{\MHM{X}}$ and we are identifying $M$ with the complex of $D^b\MHM{X}$ concentrated in degree $0$. Analogous definitions are valid for $f_*, f^*$ and $f^!$. Furthermore, these constructions imply that Theorem \ref{thm:funct-induced-functors-mhm} also holds in $K$-theory (see \cite{Schurmann:2011}) and, moreover, the natural isomorphism for cartesian squares becomes an equality.

Another important feature of mixed Hodge modules is that they actually generalize mixed Hodge structures. As mentioned in \cite{Saito:1989}, Theorem 1.4, the category of mixed Hodge modules over a single point, $\MHM{\star}$, is naturally isomorphic to the category of (rational) mixed Hodge structures $\MHS{\QQ}$. In particular, this identification endows the $K$-theory of the category of mixed Hodge modules with a natural $\K{\MHS{\QQ}} = \K{\MHM{\star}}$ module structure via the external product
$$
    \boxtimes: \K{\MHM{\star}} \times \K{\MHM{X}} \to \K{\MHM{\star \times X}} = \K{\MHM{X}}
$$
The induced functors $f_*, f_!, f^*, f^!$ commute with exterior products at the level of constructible complexes (see \cite{Schurmann:2011}), so they also commute in the category of mixed Hodge modules which means that they are $\K{\MHS{\QQ}}$-module homomorphisms. Furthermore, $f^*$ commutes with tensor products so it is also a ring homomorphism.

\begin{ex}\label{ex:cohomology-via-mhm}
Using the identification $\MHM{\star} = \MHS{\QQ}$, we can consider the Tate structure of weight 0, $\QQ_0 = \QQ(0)$, as an element of $\MHM{\star}$. By construction, this element is the unit of the ring $\K{\MHS{\QQ}}$. Moreover, for any complex algebraic variety $X$, if $c_X: X \to \star$ is the projection of $X$ onto a singleton, then the mixed Hodge module $\underline{\QQ}_X := c_X^* \QQ_0$ is the unit of $\K{\MHM{X}}$. The link between this mixed Hodge module and $X$ is that, as proven in Lemma 14.8 of \cite{Peters-Steenbrink:2008}, we can recover the compactly supported cohomology of $X$ via
$$	
	(c_X)_! \underline{\QQ}_X = [H_c^\bullet(X;\QQ)],
$$
as elements of $\K{\MHS{\QQ}}$. The analogous formula for usual cohomology and $(c_X)_*$ holds too, i.e.\ $(c_X)_* \underline{\QQ}_X = [H^\bullet(X;\QQ)]$.
\end{ex}

\begin{rmk}
Let $\pi: X \to B$ be regular fibration locally trivial in the analytic topology. In that case, the mixed Hodge module $\pi_!\underline{\QQ}_X \in \K{\MHM{B}}$ plays the role of the Hodge monodromy representation of \cite{LMN} and \cite{MM} controlling the monodromy of $\pi$.
\end{rmk}

\begin{rmk}
Let us define the semi-group homomorphism $e: \MHS{\QQ} \to \ZZ[u^{\pm 1}, v^{\pm 1}]$ that, for a mixed Hodge structure $V$ gives
$$
	\DelHod V := \sum_{p,q} \dim\left[Gr_p^{F}\left(Gr_W^{p+q}\,V\right)_\CC\right]\,u^pv^q.
$$
This map descends to a ring homomorphism $e : \K{\MHS{\QQ}} \to \ZZ[u^{\pm 1}, v^{\pm 1}]$. Now, let $X$ be a complex algebraic variety and let $[H^\bullet_c(X;\QQ)]$ be is compactly supported cohomology, as an element of $\K{\MHS{\QQ}}$. Then, we have that
$$
	\DelHod{X} = \DelHod{[H^\bullet_c(X;\QQ)]},
$$
where the left hand side is the Deligne-Hodge polynomial of $X$.
\end{rmk}

%%%%%%%%%%%%%%%%%%%%%%%%%%%%%%%%%%%%%%%%%%%%%%%%%%%%%%%%%%%%%%%%
%%%%%%%%%%%%%%%%%% SECTION: Lax TQFT %%%%%%%%%%%%%%%%%%%%%%%%%%%
%%%%%%%%%%%%%%%%%%%%%%%%%%%%%%%%%%%%%%%%%%%%%%%%%%%%%%%%%%%%%%%%

\section{Lax Topological Quantum Field Theories}

\subsection{The category of bordisms of pairs}

Let $n \geq 1$. We define the \emph{category of $n$-bordisms of pairs}, $\CBordpp{n}$ as the $2$-category (i.e. enriched category over the category of small categories, see \cite{Benabou}) given by the following data:
\begin{itemize}
	\item Objects: The objects of $\CBordpp{n}$ are pairs $(X, A)$ where $X$ is a $(n-1)$-dimensional closed oriented manifold together with a finite subset of points $A \subseteq X$ such that its intersection with each connected component of $X$ is non empty.
	\item $1$-morphisms: Given objects $(X_1, A_1)$, $(X_2, A_2)$ of $\CBordpp{n}$, a morphism $(X_1, A_1) \to (X_2, A_2)$ is a class of pairs $(W, A)$ where $W: X_1 \to X_2$ is an oriented bordism between $X_1$ and $X_2$, and $A \subseteq W$ is a finite set with $X_1 \cap A = A_1$ and $X_2 \cap A = A_2$. Two pairs $(W, A), (W',A')$ are in the same class if there exists a diffeomorphism of bordisms (i.e.\ fixing the boundaries) $F: W \to W'$ such that $F(A)=A'$.\\
With respect to the composition, given $(W, A): (X_1, A_1) \to (X_2, A_2)$ and $(W', A'): (X_2, A_2) \to (X_3, A_3)$, we define $(W',A') \circ (W,A)$ as the morphism $(W \cup_{X_2} W', A \cup A'): (X_1, A_1) \to (X_3, A_3)$ where $W \cup_{X_2} W'$ is the usual gluing of bordisms along $X_2$.
	\item $2$-morphisms: Given two $1$-morphisms $(W, A), (W',A'): (X_1, A_1) \to (X_2, A_2)$, we declare that there exists a $2$-cell $(W, A) \Rightarrow (W',A')$ if there is a diffeomorphism of bordisms $F: W \to W'$ such that $F(A) \subseteq A'$. Composition of $2$-cells is just composition of diffeomorphisms.
\end{itemize}

In this form, $\CBordpp{n}$ is not exactly a category since there is no unit morphism in the category $\Hom_{\CBordpp{n}} ((X,A),(X,A))$. This can be solved by weakening slightly the notion of bordism, allowing that $X$ itself could be seen as a bordism $X: X \to X$. With this modification, $(X,A): (X,A) \to (X,A)$ is the desired unit and it is a straightforward check to see that $\CBordpp{n}$ is a (strict) $2$-category, where strict means that the associativity axioms are satisfied ``on the nose" and not just ``up to isomorphism", in which case it is called a (weak) $2$-category.

\begin{rmk}
As a stronger version of bordisms, there is a forgetful functor $\cF: \CBordpp{n} \to \CBordp{n}$, where $\CBordp{n}$ is the usual category of oriented $n$-bordisms.
\end{rmk}

\subsection{Bimodules and soft TQFTs}
\label{subsec:bimodules}

Recall from \cite{Benabou} and \cite{Schulman} that, given a ground commutative ring $R$ with identity, we can define the $2$-category $\Bim{R}$ of $R$-algebras and bimodules whose objects are commutative $R$-algebras with unit and, given algebras $A$ and $B$, a $1$-morphism $A \to B$ is a $(A,B)$-bimodule. By convention, a $(A,B)$-bimodule is a set $M$ with a left $A$-module and a right $B$-module compatible structures, usually denoted ${_A}M_{B}$. Composition of $M: A \to B$ and $N: B \to C$ is given by ${_A}(M \otimes_B N)_C$.
	
With this definition, the set $\Hom_{\Bim{R}}(A,B)$ is naturally endowed with a category structure, namely, the category of $(A,B)$-bimodules. Hence, a $2$-morphism $M \Rightarrow N$ between $(A,B)$-bimodules is a bimodule homomorphism $f: M \to N$. Therefore, $\Bim{R}$ is a monoidal $2$-category with tensor product over $R$.

\begin{defn}
Let $R$ be a commutative ring with unit. A \emph{soft Topological Quantum Field Theory of pairs} is a lax monoidal lax functor $\sZ: \CBordpp{n} \to \Bim{R}$.
\end{defn}

Recall that a \emph{lax} functor between $2$-categories $F: \cC \to \cD$ is an assignment that:
\begin{itemize}
	\item For each object $x \in \cC$, it gives an object $F(x) \in \cD$.
	\item For each pair of objects $x, y$ of $\cC$, we have a functor
	$$
		F_{x,y}: \Hom_\cC(x, y) \to \Hom_{\cD}(F(x), F(y)).
	$$ Recall that, as $2$-category, both $\Hom_\cC(x, y)$ and $\Hom_{\cD}(F(x), F(y))$ are categories. 
	\item For each object $x \in \cC$, we have a $2$-morphism $F_{\id_x}: \id_{F(x)} \Rightarrow F_{x,x}(\id_x)$.
	\item For each triple $x,y,z \in \cC$ and every $f: x \to y$ and $g: y \to z$, we have a $2$-morphism $F_{x,y,z}(g,f): F_{y,z}(g) \circ F_{x,y}(f) \Rightarrow F_{x,z}(g \circ f)$, natural in $f$ and $g$.
\end{itemize}
Also, some technical conditions, namely the coherence conditions, have to be satisfied (see \cite{Benabou} or \cite{MacLane} for a complete definition). If the $2$-morphisms $F_{\id_x}$ and $F_{x,y,z}$ are isomorphisms, it is said that $F$ is a \emph{pseudo}-functor or a weak functor (or even simply a functor) and, if they are the identity $2$-morphism, $F$ is called a \emph{strict} functor.

Analogously, a (lax) functor $F: \cC \to \cD$ between monoidal categories is called \emph{lax monoidal} if there exists:
	\begin{itemize}
		\item A morphism $\epsilon: 1_{\cD} \to F(1_{\cC})$, where $1_{\cC}$ and $1_{\cD}$ are the units of the monoidal structure of $\cC$ and $\cD$ respectively.
		\item A natural transformation $\Delta: F(-) \otimes_{\cD} F(-) \Rightarrow F(- \otimes_\cC -)$.		
	\end{itemize}
satisfying a set of coherence conditions (for a precise definition, see Definition 1.2.10 of \cite{Leinster}). Again, if $\epsilon$ and $\Delta$ are isomorphisms, $F$ is said to be \emph{pseudo-monoidal} and, if they are identity morphisms, $F$ is called \emph{strict monoidal}, or simply monoidal.

A general recipe for building a soft TQFT from simpler data can be given as follows. Let $\cS = \Span{\CVar}$ be the (weak) $2$-category of spans of $\CVar$. As described in \cite{Benabou}, objects of $\cS$ are the same as the ones of $\CVar$. A morphism $X \to Y$ in $\cS$ between complex algebraic varieties is a pair $(f,g)$ of regular morphisms
$$
\xymatrix{
&\ar[dl]_{f} Z \ar[dr]^{g}&\\
X&& Y
}
$$
where $Z$ is a complex algebric variety. Given two morphisms $(f_1, g_1): X \to Y$ and $(f_2, g_2): Y \to Z$, say

\begin{tabular}{cc}
\begin{minipage}{3in}
$$
\xymatrix{
&\ar[dl]_{f_1} Z_1 \ar[dr]^{g_1}&\\
X&& Y
}
$$
\end{minipage}
&
\begin{minipage}{2in}
$$
\xymatrix{
&\ar[dl]_{f_2} Z_2 \ar[dr]^{g_2}&\\
Y&& Z
}
$$
\end{minipage}\vspace{0.2cm}
\end{tabular}\\
we define the composition $(f_2, g_2) \circ (f_1,g_1) = (f_1 \circ f_2', g_2 \circ g_1')$, where $f_2', g_1'$ are the morphisms in the pullback diagram
$$
\xymatrix{
&&\ar[dl]_{f'_2}W\ar[dr]^{g'_1}&&\\
&\ar[dl]_{f_1}Z_1\ar[dr]^{g_1}&&\ar[dl]_{f_2}Z_2\ar[dr]^{g_2}&\\
X&&Y&&Z
}
$$
where $W = Z_1 \times_Y Z_2$ is the product of $Z_1$ and $Z_2$ as $Y$-varieties. Finally a $2$-morphism $(f,g) \Rightarrow (f',g')$ between $X \stackrel{\hspace{3pt}f}{\leftarrow} Z \stackrel{g}{\rightarrow} Y$ and $X \stackrel{\hspace{5pt}f'}{\leftarrow} Z' \stackrel{g'}{\rightarrow} Y$ is a regular morphism $\alpha: Z' \to Z$ such that the following diagram commutes
\[
\begin{displaystyle}
   \xymatrix
   {	& Z' \ar[rd]^{g'} \ar[ld]_{f'} \ar[dd]^\alpha & \\
   		X &  & Y\\
	&Z  \ar[ru]_{g} \ar[lu]^{f}& 
   }
\end{displaystyle}   
\]
Moreover, with this definition, $\cS$ is a monoidal category with cartesian product of varieties and morphisms. 

\begin{rmk}
Actually, the same construction works verbatim for any category $\cC$ with pullbacks defining a weak $2$-category $\Span{\cC}$. Furthermore, if $\cC$ has pushouts, we can define $\Spano{\cC}$ as in the case of $\Span{\cC}$ but with all the arrows reversed. Again, if $\cC$ is a monoidal category then $\Span{\cC}$ (resp.\ $\Spano{\cC})$ also is.
\end{rmk}

Now, let $\AlgA: \CVar \to \Rng$ be a contravariant functor between the category of complex algebraic varieties and the category of (commutative unitary) rings. Set $R = \AlgA(\star)$, where $\star$ is the singleton variety, and define the lax functor $\cS \AlgA: \cS \to \Bim{R}$ as follows:
\begin{itemize}
	\item For any complex algebraic variety $X$, $\cS \AlgA(X) = \AlgA(X)$. The $R$-algebra structure on $\AlgA(X)$ is given by the morphism $R=\AlgA(\star) \to \AlgA(X)$ image of the projection $X \to \star$.
	\item Let us fix $X, Y$ algebraic varieties. We define 
	$$\cS \AlgA_{X,Y}: \Hom_{\cS}(X,Y) \to \Hom_{\Bim{R}}(\AlgA(X), \AlgA(Y))$$ as the covariant functor that:
		\begin{itemize}
			\item For any $1$-morphism $X \stackrel{f}{\leftarrow} Z \stackrel{g}{\rightarrow} Y$ it assigns $\cS \AlgA_{X,Y}(f,g) = \AlgA(Z)$ with the $(\AlgA(X),  \AlgA(Y))$-bimodule structure given by $az = \AlgA(f)(a) \cdot z$ and $zb = z \cdot \AlgA(g)(b)$ for $a \in \AlgA(X)$, $b \in \AlgA(Y)$ and $z \in \AlgA(Z)$.
			\item For a $2$-morphism $\alpha: (f,g) \Rightarrow (f',g')$ given by a regular morphism $\alpha: Z' \to Z$, we define $\cS \AlgA_{X,Y}(\alpha) = \AlgA(\alpha): \AlgA(Z) \to \AlgA(Z')$. Since $\AlgA$ is a functor to rings, $\AlgA(\alpha)$ is a bimodule homomorphism.
		\end{itemize}
	\item The $2$-morphism $\cS \AlgA_{\id_X}: \id_{\cS \AlgA(X)} \Rightarrow \cS \AlgA_{X,X}(\id_X)$ is the identity.
	\item Let $X \stackrel{f_1}{\leftarrow} Z_1 \stackrel{g_1}{\rightarrow} Y$ and $Y \stackrel{f_2}{\leftarrow} Z_2 \stackrel{g_2}{\rightarrow} Z$. By definition of $\cS$,
\begin{align*}
\cS \AlgA_{Y, Z}(f_2,g_2) \circ \cS \AlgA_{Y, X}(f_1, g_1) &= {_{\AlgA(Y)}}\AlgA(Z_2)_{\AlgA(Z)} \circ {_{\AlgA(X)}}A(Z_1)_{\AlgA(Y)} \\
	&= {_{\AlgA(X)}}\left[\AlgA(Z_1) \otimes_{\AlgA(Y)} \AlgA(Z_2)\right]_{\AlgA(Z)} \\
\cS \AlgA_{X,Z}((f_2,g_2) \circ (f_1, g_1)) &= {_{\AlgA(X)}}\left[\AlgA(Z_1 \times_Y Z_2)\right]_{\AlgA(Z)}
\end{align*}
Hence, the $2$-morphism $\cS \AlgA_{X,Y,Z}((f_2, g_2), (f_1, g_1))$ is the $(\AlgA(X), \AlgA(Z))$-bimodule homomorphism
$$
	\AlgA(Z_1) \otimes_{\AlgA(Y)} \AlgA(Z_2) \to \AlgA(Z_1 \times_Y Z_2)
$$
given by $z_1 \otimes z_2 \mapsto \AlgA(p_1)(z_1) \cdot \AlgA(p_2)(z_2)$ where $p_1$ and $p_2$ are the projections $p_1: Z_1 \times_Y Z_2 \to Z_1$, $p_2: Z_1 \times_Y Z_2 \to Z_2$.
\end{itemize}

Furthermore, with the construction above, the lax functor $\cS \AlgA: \cS \to \Bim{R}$ is also lax monoidal. Since $\cS \AlgA$ preserves the units of the monoidal structures it is enough to define the morphisms $\Delta_{X,Y}: \AlgA(X) \otimes_R \AlgA(Y) \to \AlgA(X \times Y)$. For that, just take as bimodule the own ring $\AlgA(X \times Y)$, where the left $(\AlgA(X) \otimes_R \AlgA(Y))$-module structure commes from the external product $\AlgA(X) \otimes_R \AlgA(Y) \to \AlgA(X \times Y)$, $z_1 \otimes z_2 \mapsto \AlgA(p_1)(z_1) \cdot \AlgA(p_2)(z_2)$ (being $p_1, p_2$ the respective projections) which is a ring homomorphism.

\begin{rmk}
If the functor $\AlgA$ is monoidal, then $\cS \AlgA$ is strict monoidal. Moreover, if it satisfies that $\AlgA(Z_1) \otimes_{\AlgA(Y)} \AlgA(Z_2) = \AlgA(Z_1 \times_Y Z_2)$ (resp.\ isomorphic) then it is also a strict (resp.\ pseudo) functor.
\end{rmk}

Therefore, with this construction at hand, given a lax monoidal lax functor
$$
	\Geo: \CBordpp{n} \to \cS,
$$
called the \emph{geometrisation}, and a functor
$$
	\AlgA: \CVar \to \Rng,
$$
called the \emph{algebraisation}, we can build a soft TQFT, $\sZ = \sZ_{\Geo,\AlgA}$, by
$$
	\sZ_{\Geo,\AlgA} = \cS \AlgA \circ \Geo: \CBordpp{n} \to \Bim{R}.
$$
We will use this approach in section \ref{sec:LTQFT-E-pol} to define a soft TQFT generalizing Deligne-Hodge polynomials of representation varieties.

\subsection{$2$-categories of modules with twists}
\label{subsec:modules-twists}

Let $R$ be a fixed ring (commutative and with unit). Given two homomorphisms of $R$-modules $f,g: M \to N$, we say that $g$ is an \emph{immediate twist} of $f$ if there exists an $R$-module $L$, homomorphisms $f_1: M \to L$, $f_2: L \to N$ and $\psi: D \to D$ such that $f = f_2 \circ f_1$ and $g = f_2 \circ \psi \circ f_1$.
\[
\begin{displaystyle}
   \xymatrix
   {	M \ar[r]^{f_1} \ar@/_1pc/[rr]_g & L \ar[r]^{f_2} \ar@(ul,ur)^{\psi} & N
   }
\end{displaystyle}   
\]
In general, given $f,g: M \to N$ two $R$-module homomorphisms, we say that $g$ is a \emph{twist} of $f$ if there exists a finite sequence $f=f_0, f_1, \ldots, f_r = g: M \to N$ of homomorphisms such that $f_{i+1}$ is an immediate twist of $f_i$.

\begin{defn}
\label{defn:twisted-modules}
Let $R$ be a ring. The $2$-category of $R$-modules with twists, $\Mod{R}_t$ is the category whose objects are $R$-modules, its $1$-morphisms are $R$-modules homomorphisms and, given homomorphisms $f$ and $g$, we have a $2$-morphism $f \Rightarrow g$ if and only if $g$ is a twist of $f$. Moreover, $\Mod{R}_t$ is a monoidal category with the usual tensor product.
\end{defn}

\begin{defn}Let $R$ be a commutative ring with unit.
A \emph{lax monoidal Topological Quantum Field Theory of pairs} is a lax monoidal strict $2$-functor $\lZ: \CBordpp{n} \to \Mod{R}_t$.
\end{defn}

\begin{rmk}
In the literature, it is customary to forget about the $2$-category structure and to say that a lax monoidal TQFT is just a lax monoidal functor between $\CBordp{n}$ and $\Mod{R}$. In this paper, the $2$-category structures are chosen to suit the geometric situation.
\end{rmk}

\begin{rmk}
Since we requiere $\lZ$ to be only lax monoidal, some of the properties of Topological Quantum Field Theories can be lost. For example, duality arguments no longer hold and in particular, $\lZ(X)$ can be not finitely generated. This is the case of the construction of section \ref{sec:LTQFT-E-pol}.
\end{rmk}

A lax monoidal TQFT is, in some sense, stronger than a soft TQFT. As described in section \ref{subsec:bimodules} suppose that we have a contravariant functor $\AlgA: \CVar \to \Rng$ and set $R=\AlgA(\star)$.
Furthermore, suppose that, together with this functor, we have a covariant functor $\AlgB: \CVar \to \Mod{R}$ (being $\Mod{R}$ the usual category of $R$-modules) such that $\AlgB(X)=\AlgA(X)$ for any complex algebraic variety $X$ and satisfying the \emph{Beck-Chevalley condition}, sometimes called the \emph{base change formula}, which requires that, for any pullback diagram in $\CVar$

\[
\begin{displaystyle}
   \xymatrix
   {
   		Y \ar[d]_{f'} \ar[r]^{g'}& Y_1 \ar[d]^{f} \\
   		Y_2 \ar[r]_{g} & X
   }
\end{displaystyle}   
\]
we have $\AlgA(g)\AlgB(f) = \AlgB(f')\AlgA(g')$.

In this framework, we can define the strict $2$-functor $\cS_{\AlgA,\AlgB}: \cS=\Span{\CVar} \to \Mod{R}_t$ as follows:
\begin{itemize}
	\item For any algebraic variety $X$ we define $\cS_{\AlgA,\AlgB}(X)=\AlgA(X)$.
	\item Fixed $X, Y$ complex algebraic varieties, we define the functor
	$$
		(\cS_{\AlgA,\AlgB})_{X,Y}: \Hom_{\cS}(X, Y) \to \Hom_{\Mod{R}_t}(\AlgA(X), \AlgA(Y))$$ by:
	\begin{itemize}
		\item For a $1$-morphism $X \stackrel{f}{\leftarrow} Z \stackrel{g}{\rightarrow} Y$ we define $\cS_{\AlgA,\AlgB}(f,g) = \AlgB(g) \circ \AlgA(f): \AlgA(X) \to \AlgA(Z) \to \AlgA(Y)$.
		\item For a $2$-morphism $\alpha: (f,g) \Rightarrow (f',g')$ given by $\alpha: Z' \to Z$ we define the immediate twist $\psi = \AlgB(\alpha) \circ \AlgA(\alpha): \AlgA(Z) \to \AlgA(Z)$. Since $\alpha$ is a $2$-cell in $\cS$ we have that $\AlgB(g)\,\psi\, \AlgA(f) = \AlgB(g') \AlgA(f')$. Observe that, if $\alpha$ is an isomorphism, then $\psi = \id$.
	\end{itemize}
	\item For $(\cS_{\AlgA,\AlgB})_{\id_X}$ we take the identity $2$-cell.
	\item Given $1$-morphisms $(f_1, g_1): X \to Y$ and $(f_2, g_2): Y \to Z$ we have $(\cS_{\AlgA,\AlgB})_{Y,Z}(f_2,g_2) \circ (\cS_{\AlgA,\AlgB})_{X,Y}(f_1,g_1) = \AlgB(g_2)\AlgA(f_2)\AlgB(g_1)\AlgA(f_1)$. On the other hand, $(\cS_{\AlgA,\AlgB})_{X,Z}((f_2,g_2) \circ (f_1, g_1))  =  \AlgB(g_2) \AlgB(g_1')  \AlgA(f_2')  \AlgA(f_1)$, where $g_1$ and $f_2'$ are the maps in the pullback
	$$
\xymatrix{
&&\ar[dl]_{f'_2}Z_1\times_Y Z_2\ar[dr]^{g'_1}&&\\
&\ar[dl]_{f_1}Z_1\ar[dr]^{g_1}&&\ar[dl]_{f_2}Z_2\ar[dr]^{g_2}&\\
X&&Y&&Z.
}
$$
By the Beck-Chevalley condition we have $\AlgB(g_1')\AlgA(f_2') = \AlgA(f_2)\AlgB(g_1)$ and the two morphisms agree. Thus, we can take the $2$-cell $(\cS_{\AlgA,\AlgB})_{X,Y,Z}$ as the identity.
\end{itemize}

\begin{rmk}
The Beck-Chevalley condition appears naturally in the context of Grothendieck's yoga of six functors $f_*, f^*, f_!, f^!, \otimes$ and $\mathbb{D}$ in which $(f_*, f^*)$ and $(f_!, f^!)$ are adjoints, and $f^*$ and $f_!$ satisfy the Beck-Chevalley condition. In this context, we can take $\AlgA$ to be the functor $f \mapsto f^*$ and $\AlgB$ the functor $f \mapsto f_!$ as we will do in section \ref{sec:LTQFT-E-pol}. For further infomation, see for example \cite{Fausk-May:2003} or \cite{Ayoub}.
\end{rmk}

\begin{rmk}
If the Beck-Chevalley condition was satisfied up to natural isomorphism, then this corresponds to equality after a pair of automorphisms in $\AlgA(Z_1)$ and $\AlgA(Z_2)$ which corresponds to an invertible $2$-cell in $\Mod{R}_t$. In this case, $\cS_{\AlgA,\AlgB}$ would be a pseudo-functor. We will not need this trick in this paper.
\end{rmk}

As in the previous case, $\cS_{\AlgA,\AlgB}: \cS \to \Mod{R}_t$ is automatically lax monoidal taking $\Delta_{X, Y}: \cS_{\AlgA,B} (X) \otimes \cS_{\AlgA,\AlgB}(Y) \to \cS_{\AlgA,\AlgB}(X \times Y)$ to be the external product with respect to $\AlgA$.
Therefore, given a geometrisation functor, i.e.\ a lax monoidal strict functor $\Geo: \CBordpp{n} \to \cS$, and these two functors $\AlgA: \CVar \to \Rng$ and $\AlgB: \CVar \to \Mod{R}$ as algebraisations, we can build a lax monoidal TQFT, $\lZ=\lZ_{\Geo, \AlgA, \AlgB}$ by
$$
	\lZ_{\Geo, \AlgA, \AlgB} = \cS_{\AlgA,\AlgB} \circ \Geo: \CBordpp{n} \to \Mod{R}_t.
$$
Observe that the previous soft TQFT, $\sZ = \cS \AlgA \circ \Geo$ can also be constructed in this setting, so in this sense $\lZ$ is stronger than $\sZ$.

\subsection{The parabolic case}
\label{subsec:parabolic-case}
For some applications, it is useful to consider the so-called parabolic case.
In this setting, the starting point is a fixed set $\Lambda$ that we will call the \emph{parabolic data}. Given a compact $n$-dimensional manifold $W$, possibly with boundary, we will denote by $\Par{W}$ the set of closed connected subvarieties $S \subseteq W$ of dimension $n-2$ such that $S \cap \partial W = \emptyset$. A \emph{parabolic structure}, $Q$, on $W$ is a finite set (possibly empty) of pairs $Q = \left\{(S_1,\lambda_1), \ldots, (S_s, \lambda_s)\right\}$ where $S_i \in \Par{W}$ and $\lambda_i \in \Lambda$.

With this notion, we can improve our previous category $\CBordpp{n}$ to the $2$-category of pairs of bordisms with parabolic data $\Lambda$, $\CBordpp{n}(\Lambda)$. The objects of this category are the same than for $\CBordpp{n}$. For morphisms, given objects $(X_1, A_1)$ and $(X_2, A_2)$ of $\CBordpp{n}$, a morphism $(X_1, A_1) \to (X_2, A_2)$ is a class of triples $(W, A, Q)$ where $(W, A): (X_1, A_1) \to (X_2, A_2)$ is a bordism of pairs and $Q$ is a parabolic structure on $W$ such that $S \cap A = \emptyset$ for any $(S, \lambda) \in Q$. As in the case of $\CBordpp{n}$, two triples $(W, A, Q)$ and $(W', A', Q')$ are in the same class if there exists a diffeomorphism $F: W \to W'$ preserving the boundaries such that $F(A)=A'$ and $(S, \lambda) \in Q$ if and only if $(F(S), \lambda) \in Q'$. As expected, composition of morphisms in $\CBordpp{n}(\Lambda)$ is defined by $(W', A', Q') \circ (W, A, Q) = (W \cup W', A \cup A', Q \cup Q')$. In the same spirit, we have a $2$-morphism $(W, A, Q) \Rightarrow (W', A', Q')$ if there exists a boundary preserving diffeomorphism $F: W \to W'$ such that $F(A) \subseteq A'$ and $(S, \lambda) \in Q$ if and only if $(F(S), \lambda) \in Q'$.

Then, in analogy with the non-parabolic case, a lax monoidal lax functor $\sZ: \CBordpp{n}(\Lambda) \to \Bim{R}$ is called a parabolic soft TQFT and a lax monoidal strict functor $\lZ: \CBordpp{n}(\Lambda) \to \Mod{R}_t$ is called a parabolic lax monoidal TQFT.

%%%%%%%%%%%%%%%%%%%%%%%%%%%%%%%%%%%%%%%%%%%%%%%%%%%%%%%%%%%%%%%%
%%%%%%%%%%%%%%%%%% SECTION: LTQFT FOR E-POL %%%%%%%%%%%%%%%%%%%%%
%%%%%%%%%%%%%%%%%%%%%%%%%%%%%%%%%%%%%%%%%%%%%%%%%%%%%%%%%%%%%%%%

\section{TQFTs for Deligne-Hodge polynomials of representations varieties}
\label{sec:LTQFT-E-pol}

In this section, we will use the previous ideas to construct a soft TQFT and a lax monoidal TQFT allowing computation of Deligne-Hodge polynomials of representation varieties. Both theories share a common geometrisation functor, constructed via the fundamental groupoid functor as described in section \ref{subsec:E-pol-geometric-part}. On the other hand, for the algebraisation we will use, in section \ref{subsec:E-pol-algebraic-part}, the properties of mixed Hodge modules so that $\AlgA$ will be the pullback of mixed Hodge modules and $\AlgB$ the pushforward.

\subsection{The geometrisation functor}
\label{subsec:E-pol-geometric-part}

Recall from \cite{Brown}, Chapter 6, that given a topological space $X$ and a subset $A \subseteq X$, we can define the \emph{fundamental groupoid of $X$ with respect to $A$}, $\Pi(X, A)$, as the category whose objects are the points in $A$ and, given $a, b \in A$, $\Hom(a,b)$ is the set of homotopy classes (with fixed endpoints) of paths between $a$ and $b$. It is a straighforward check that this category is actually a groupoid and only depends on the homotopy type of the pair $(X,A)$. In particular, if $A = \left\{x_0\right\}$, $\Pi(X, A)$ has a single object whose automorphism group is $\pi_1(X, x_0)$, the fundamental group of $X$ based on $x_0$. For convenience, if $A$ is any set, not necessarily a subset of $X$, we will denote $\Pi(X, A) := \Pi(X, X \cap A)$ and we declare that $\Pi(\emptyset, \emptyset)$ is the singleton category.

With this notion at hand, we define the strict $2$-functor $\Pi: \CBordpp{n} \to \Spano{\Grpd}$, where $\Grpd$ is the category of groupoids, as follows:
\begin{itemize}
	\item For any object $(X, A)$ of $\CBordpp{n}$ we set $\Pi(X,A)$ as the fundamental groupoid of $(X,A)$.
	\item Given objects $(X_1, A_1), (X_2, A_2)$, the functor
	$$\Pi_{(X_1 A_1), (X_2, A_2)}: \Hom_{\CBordpp{n}}((X_1, A_1), (X_2, A_2)) \to \Hom_{\Spano{\Grpd}}(\Pi(X_1,A_1), \Pi(X_2,A_2))
	$$ is given by:
	\begin{itemize}
		\item For any $1$-morphism $(W,A): (X_1, A_1) \to (X_2, A_2)$ we define $\Pi_{(X_1 A_1), (X_2, A_2)}(W,A)$ to be the span
		$$
			\Pi(X_2,A_2) \stackrel{i_2}{\longrightarrow} \Pi(W,A) \stackrel{i_1}{\longleftarrow} \Pi(X_1,A_1)
		$$
where $i_1, i_2$ are the induced functions on the level of groupoids by the inclusions of pairs $(X_1, A_1), (X_2, A_2) \hookrightarrow (W,A)$.
		\item For a $2$-morphism $(W, A) \Rightarrow (W', A')$ given by a diffeomorphism $F: W \to W'$, we obtain a groupoid homomorphism $\Pi F: \Pi(W, A) \to \Pi(W',A')$ giving rise to a commutative diagram
\[
\begin{displaystyle}
   \xymatrix
   {	& \Pi(W, A)   \ar[d]^{\Pi F} & \\
   		\Pi(X_2, A_2) \ar[ru]^{i_2} \ar[r]_{i_2'} & \Pi(W', A')   & \Pi(X_1, A_1) \ar[lu]_{i_1}  \ar[l]^{i_1'}
   }
\end{displaystyle}   
\]
which is a $2$-cell in $\Spano{\Grpd}$.
	\end{itemize}
	\item For $\Pi_{\id_{(X,A)}}$ we take the identity.
	\item With respect to composition, let us take $(W, A): (X_1, A_1) \to (X_2, A_2)$ and $(W', A'): (X_2, A_2) \to (X_3, A_3)$ two $1$-morphisms. Set $W'' = W \cup_{X_2} W'$ and $A'' = A \cup A'$ so that $(W', A') \circ (W, A) = (W'', A'')$. In order to identify $\Pi(W'', A'')$, let $V \subseteq W''$ be an open bicollar of $X_2$ such that $V \cap A'' = A_2$. Set $U_1 = W \cup V$ and $U_2 = W' \cup V$.
	
	By construction, $\left\{U_1, U_2\right\}$ is an open covering of $W''$ such that $(U_1, A'' \cap U_1)$ is homotopically equivalent to $(W, A)$, $(U_2, A'' \cap U_2)$ is homotopically equivalent to $(W', A')$ and $(U_1 \cap U_2, A'' \cap U_1 \cap U_2)$ is homotopically equivalent to $(V, A'' \cap V)$ which is homotopically equivalent to $(X_2, A_2)$. Therefore, by Seifert-van Kampen theorem for fundamental groupoids (see \cite{Brown}, \cite{Brown:1967} and \cite{Higgins}) we have a pushout diagram induced by inclusions
	\[
\begin{displaystyle}
   \xymatrix
   {	\Pi(U_1 \cap U_2, A'')=\Pi(X_2, A_2) \ar[r] \ar[d] & \Pi(U_1, A'') = \Pi(W, A) \ar[d] \\
   		\Pi(U_2, A'') = \Pi(W', A') \ar[r]& \Pi(W'', A'')
   }
\end{displaystyle}   
\]
But observe that, by definition, $\Pi(W', A') \circ \Pi(W, A)$ is precisely this pushout, so we can take the functors $\Pi_{X_1,X_2,X_3}$ as the identities.
\end{itemize}

We can slightly improve the previous construction. Given a groupoid $\cG$, we will say that $\cG$ is \emph{finitely generated} if $\Obj{\cG}$ is finite and, for any object $a$ of $\cG$, $\Hom_{\cG}(a,a)$ (usually denoted $\cG_a$, the vertex group of $a$) is a finitely generated group. We will denote by $\Grpdo$ the category of finitely generated groupoids.

\begin{rmk}
Given a groupoid $\cG$, two objects $a,b \in \cG$ are said to be connected if $\Hom_{\cG}(a,b)$ is not empty. In particular, this means that $\cG_a$ and $\cG_b$ are isomorphic groups so, in order to check whether $\cG$ is finitely generated, it is enough to check it on a point of every connected component.
\end{rmk}

\begin{rmk}
Let $X$ be a compact connected manifold (possibly with boundary) and let $A \subseteq X$ be finite. As we mentioned above, for any $a \in A$, $\Pi(X,A)_a = \pi_1(X, a)$. But a compact connected manifold has the homotopy type of a finite CW-complex, so in particular has finitely generated fundamental group. Hence, $\Pi(X, A)$ is finitely generated. Therefore, the previous functor actually can be promoted to a functor $\Pi: \CBordpp{n} \to \Spano{\Grpdo}$.
\end{rmk}

Now, let $G$ be a complex algebraic group. Seeing $G$ as a groupoid, we can consider the functor $\Hom_{\Grpd}(-, G): \Grpd \to \Set$. Moreover, if $\cG$ is finitely generated, then $\Hom(\cG, G)$ has a natural structure of complex algebraic variety. To see that, pick a set $S = \left\{a_1, \ldots, a_s\right\}$ of objects of $\cG$ such that every connected component contains exactly one element of $S$. Moreover, for any object $a$ of $\cG$, pick a morphism $f_a: a \to a_i$ where $a_i$ is the object of $S$ in the connected component of $a$. Hence, if $\rho: \cG \to G$ is a groupoid homomorphism, it is uniquely determined by the group representations $\rho_i: \cG_{a_i} \to G$ for $a_i \in S$ together with elements $g_a$ corresponding to the morphisms $f_a$ for any object $a$. Since the elements $g_a$ can be chosen without any restriction, if $\cG$ has $n$ objects, we have
$$
	\Hom(\cG, G) \cong \Hom(\cG_{a_1},G) \times \ldots \times \Hom(\cG_{a_1}, G) \times G^{n-s}
$$
and each of these factors has a natural structure of complex algebraic variety as representation variety. This endows $\Hom(\cG, G)$ with an algebraic structure which can be shown not to depend on the choices.

\begin{defn}
Let $(W, A)$ be a pair of topological spaces such that $\Pi(W,A)$ is finitely generated (for example, if $W$ is a compact manifold and $A$ is finite). We denote the variety $\Rep{G}(W,A) = \Hom(\Pi(W, A), G)$ and we call it the \emph{representation variety} of $(W, A)$ into $G$. In particular, if $A$ is a singleton, we recover the usual representation varieties of group homomorphisms $\rho: \pi_1(W) \to G$, just denoted by $\Rep{G}(W)$.
\end{defn}

\begin{rmk}
If we drop out the requirement of $\cG$ being finitely generated, we can still endow $\Hom(\cG, G)$ with a scheme structure following the same lines (see, for example \cite{Nakamoto}). However, in general, this scheme is no longer of finite type. For this reason, in the definition of $\CBordpp{n}$, we demand the subset $A \subseteq W$ to be finite.
\end{rmk}

Therefore, we can promote this functor to a contravariant functor $\Hom(-,G): \Grpdo \to \CVar$. Recall that $\Hom(-,G)$ sends colimits of $\Grpd$ into limits of $\CVar$ so, in particular, sends pushouts into pullback and, thus, defines a functor $\Hom(-,G): \Spano{\Grpdo} \to \Span{\CVar}$. With this functor at hand, we can finally define the geometrisation functor as
$$
	\Geo = \Hom(-,G) \circ \Pi: \CBordpp{n} \to \cS = \Span{\CVar}.
$$
Observe that, since $\Pi$ and $\Hom(-,G)$ are both (strict) monoidal functors then $\Geo$ is strict monoidal.

\subsection{The algebraisation functor}
\label{subsec:E-pol-algebraic-part}

As algebraisation functor for our soft TQFT, let us consider the contravariant functor $\AlgA: \CVar \to \Rng$ that, for a complex algebraic variety $X$ gives $\AlgA(X)=\K{\MHM{X}}$ and, for a regular morphism $f: X \to Y$ it assigns $\AlgA(f) = f^*: \K{\MHM{Y}} \to \K{\MHM{X}}$. As we mention in section \ref{subsec:mixed-hodge-modules}, $f^*$ is a ring homomorphism and $\AlgA(\star) = \K{\MHM{\star}} = \K{\MHS{\QQ}}$. With these choices, the corresponding soft TQFT is $\sZ = \cS \AlgA \circ \Geo: \CBordpp{n} \to \Bim{R}$ with $R = \K{\MHS{\QQ}}$. As described in section \ref{subsec:bimodules}, it satisfies:
\begin{itemize}
	\item For any pair $(X, A)$, where $X$ is a compact $(n-1)$-dimensional manifold and $A \subseteq X$ is finite, we have $\sZ(X,A) = \K{\MHM{\Rep{G}(X,A)}}$.
	\item For a $1$-morphism $(W, A): (X_1, A_1) \to (X_2, A_2)$ we assign
	$$
		\sZ(W, A) = \K{\MHM{\Rep{G}(W,A)}}
	$$
with the structure of a $(\K{\MHM{\Rep{G}(X_1,A_1)}}, \K{\MHM{\Rep{G}(X_2,A_2)}})$-bimodule. In particular, taking the unit $\underline{\QQ} \in \K{\MHM{\Rep{G}(W,A)}}$ and the projection onto a singleton $c: {\Rep{G}(W,A)} \to \star$, by Example \ref{ex:cohomology-via-mhm} we have
$$
	{c}_! \underline{\QQ} = \left[H^\bullet_c({\Rep{G}(W,A)}; \QQ)\right]
$$
as mixed Hodge structures.
	\item For a $2$-morphism $\alpha: (W, A) \Rightarrow (W', A')$ we obtain a bimodule homomorphism
	$$
	\sZ(\alpha): \K{\MHM{\Rep{G}(W',A')}} \to \K{\MHM{\Rep{G}(W,A)}}
	$$
\end{itemize}

\begin{rmk}
As we mentioned in section \ref{subsec:bimodules}, $\cS \AlgA: \Span{\CVar} \to \Bim{R}$ is automatically a lax monoidal lax functor. However, it is not strict monoidal since, in general, $\K{\MHM{X \times Y}} \neq \K{\MHM{X}} \otimes_{R} \K{\MHM{Y}}$.
\end{rmk}

Furthermore, we can also consider the covariant functor $\AlgB: \CVar \to \Mod{R}$, where $R = \AlgA(\star)=\K{\MHS{\QQ}}$ given by $\AlgB(X) = \K{\MHM{X}}$ for an algebraic variety $X$ and, for a regular morphism $f: X \to Y$, it assigns $\AlgB(f) = f_!: \K{\MHM{X}} \to \K{\MHM{Y}}$. In this case, as described in section \ref{subsec:modules-twists}, the corresponding lax monoidal TQFT is $\lZ = \cS_{\AlgA,\AlgB} \circ \Geo: \CBordpp{n} \to \Mod{R}_t$. It assigns:
\begin{itemize}
	\item For a pair $(X, A)$, $\lZ(X,A) = \K{\MHM{\Rep{G}(X,A)}}$.
	\item For a $1$-morphism $(W, A): (X_1, A_1) \to (X_2, A_2)$, it assigns the $R$-module homomorphism
	$$
		\lZ(W, A): \K{\MHM{\Rep{G}(X_1,A_1)}} \stackrel{i_1^*}{\longrightarrow} \K{\MHM{\Rep{G}(W,A)}} \stackrel{{i_2}_!}{\longrightarrow} \K{\MHM{\Rep{G}(X_2,A_2)}}
	$$
	\item The existence of a $2$-morphism $\alpha: (W, A) \Rightarrow (W', A')$ implies that $\lZ(W', A') = ({i_2'})_!\circ (i_1')^*: \K{\MHM{\Rep{G}(X_1,A_1)}} \to \K{\MHM{\Rep{G}(X_2,A_2)}}$ can be obtained from $\lZ(W, A) = ({i_2})_!\circ (i_1)^*: \K{\MHM{\Rep{G}(X_1,A_1)}} \to \K{\MHM{\Rep{G}(X_2,A_2)}}$ by twists.
\end{itemize}

In particular, let $W$ be a connected closed oriented $n$-dimensional manifold and let $A \subseteq W$ finite. We can see $W$ as a $1$-morphism $(W, A): (\emptyset, \emptyset) \to (\emptyset, \emptyset)$. In that case, since $\lZ(\emptyset, \emptyset) = \K{\MHM{\Rep{G}(\emptyset, \emptyset)}} = \K{\MHM{\star}} = \K{\MHS{\QQ}}$ we obtain that $\lZ(W, A)$ is the morphism
$$
	\lZ(W, A): \K{\MHS{\QQ}} \stackrel{c^*}{\longrightarrow} \K{\MHM{\Rep{G}(W, A)}} \stackrel{c_!}{\longrightarrow} \K{\MHS{\QQ}},
$$
where $c: \Rep{G}(W, A) \to \star$ is the projection onto a point. In particular, for the unit $\QQ_0 \in \MHS{\QQ}$ we have
$$
	\lZ(W, A)(\QQ_0) = c_!\,c^*\QQ_0 = c_!\,\underline{\QQ}_{\Rep{G}(W, A)} = \left[H^\bullet_c\left(\Rep{G}(W, A) ;\QQ\right)\right].
$$
Hence, taking into account that $\Rep{G}(W, A) = \Rep{G}(W) \times G^{|A|-1}$ we have that
$$
	\lZ(W, A)(\QQ_0) = \left[H_c^\bullet\left(\Rep{G}(W); \QQ\right)\right] \otimes \left[H_c^\bullet\left(G; \QQ\right)\right]^{|A|-1}
$$
Therefore, we have proved the main result of this paper in the non-parabolic case.

\begin{thm}\label{thm:existence-s-LTQFT}
There exists a lax monoidal TQFT, $\lZ: \CBordpp{n} \to \Mod{\K{\MHS{\QQ}}}_t$ such that, for any $n$-dimensional connected closed orientable manifold $W$ and any non-empty finite subset $A \subseteq W$ we have
$$
	\lZ(W, A)(\QQ_0) = \left[H_c^\bullet\left(\Rep{G}(W); \QQ\right)\right] \otimes \left[H_c^\bullet\left(G; \QQ\right)\right]^{|A|-1}
$$ 
where $\QQ_0 \in \K{\MHS{\QQ}}$ is the unit Hodge structure. In particular, this means that
$$
	\DelHod{\lZ(W, A)(\QQ_0)} = \DelHod{G}^{|A|-1} \DelHod{\Rep{G}(W)}.
$$ 
\end{thm}

\begin{rmk}
An analogous formula holds in the non-connected case. Suppose that $W = W_1 \sqcup \ldots \sqcup W_s$ with $W_i$ connected and denote $A_i = A \cap W_i$. Then $\Rep{G}(W, A) = \Rep{G}(W_1, A_1) \times \ldots \times \Rep{G}(W_s, A_s)$ so $\lZ(W, A) = \bigotimes_{i} \lZ(W_i, A_i)$ and, thus
$$
	\DelHod{\lZ(W, A)(\QQ_0)} = \DelHod{G}^{|A|-s} \prod_{i=1}^s\DelHod{\Rep{G}(W_i)}.
$$
\end{rmk}

\subsection{Almost-$\textrm{TQFT}$ and computational methods}
\label{subsec:almost-TQFT}

For computational purposes, we can restrict our attention to a wide subcategory of bordisms. Along this section, we will forget about the $2$-category structure of $\CBordpp{n}$ and we will see it as a $1$-category.

First of all, let us consider the subcategory $\CTubppo{n} \subseteq \CBordpp{n}$ of \emph{strict tubes of pairs}. An object $(X,A)$ of $\CBordpp{n}$ is an object of $\CTubppo{n}$ if $X$ is connected or empty. Given objects $(X_1, A_1)$ and $(X_2, A_2)$ of $\CTubppo{n}$, a morphism $(W, A): (X_1, A_1) \to (X_2, A_2)$ of $\CBordpp{n}$ is in $\CTubppo{n}$ if $W$ is connected. We will call such morphisms strict tubes.

From this category, we consider the \emph{category of tubes}, $\CTubpp{n}$ as the subcategory of $\CBordpp{n}$ whose objects and morphisms are disjoint unions of the ones of $\CTubppo{n}$. Observe that, in particular, $\CTubpp{n}$ is a wide subcategory of $\CBordpp{n}$ i.e.\ they have the same objects. As in the case of $\CBordpp{n}$, $\CTubpp{n}$ is a monoidal category with the disjoint union.

\begin{defn}
Let $R$ be a ring. An \emph{almost-TQFT of pairs} is a monoidal functor $\mathfrak{Z}: \CTubpp{n} \to \Mod{R}$.
\end{defn}

\begin{rmk}
Since $\CTubpp{n}$ does not contain all the bordisms, dualizing arguments no longer hold for almost-TQFTs. For example, for $n=2$, the pair of pants is not a tube, so, in contrast with TQFTs, we cannot assure that $\mathfrak{Z}(S^1, \star)$ is a Frobenius algebra. In particular, $\mathfrak{Z}(S^1, \star)$ is not forced to be finite dimensional.
\end{rmk}

\begin{rmk}
An almost-TQFT gives an effective way of computing invariants as follows. Fix $n \geq 1$ and suppose that we have a set of generators $\Delta$ for the morphisms of $\CTubppo{n}$, i.e.\ after boundary orientation preserving diffeomorphisms, every morphism of $\CTubppo{n}$ is a compositions of elements of $\Delta$. These generators can be obtained, for example, by means of Morse theory (see Section 1.4 of \cite{Kock:2004} for $n=2$).

Suppose that we want to compute an invariant that, for a closed connected orientable $n$-dimensional manifold $W$ and a finite set $A \subseteq W$ is given by $\mathfrak{Z}(W,A)(1)$. In that case, seeing $(W,A)$ as a morphism $(W,A): \emptyset \to \emptyset$, we can decompose $(W,A) = W_{s} \circ \ldots \circ W_{1}$ with $W_{i} \in \Delta$. Thus, for $\mathfrak{Z}(W, A): R \to R$ we have
$$
	\mathfrak{Z}(W, A)(1) = \mathfrak{Z}(W_{s}) \circ \ldots \circ \mathfrak{Z}(W_{1})(1).
$$
Hence, the knowledge of $\mathfrak{Z}(W_i)$ for $W_i \in \Delta$ is enough to compute that invariant for closed manifolds.
\end{rmk}

Given a lax monoidal TQFT, $\lZ: \CBordpp{n} \to \Mod{R}_t$, we can build a natural almost-TQFT from it, $\mathfrak{Z} = \mathfrak{Z}_{\lZ}: \CTubpp{n} \to \Mod{R}$. For that, forgetting about the $2$-category structure, the restriction of $\lZ$ to $\CTubppo{n}$ gives us a functor $\lZ: \CTubppo{n} \to \Mod{R}$. Now, we define $\mathfrak{Z}$ by:
\begin{itemize}
	\item Given objects $(X^i, A^i)$ of $\CTubppo{n}$, we take
	$$
		\mathfrak{Z}\left(\bigsqcup_i \,(X^i, A^i)\right) = {\bigotimes_i}\, \lZ(X^i, A^i)
	$$
where the tensor product is taken over $R$.
	\item Given strict tubes $(W^i, A^i): (X_{1}^i, A_1^i) \to (X_2^i, A_2^i)$, we define
	$$
		\mathfrak{Z}\left(\bigsqcup_i \,(W^i, A^i)\right) = {\bigotimes_i}\, \lZ(W^j, A^j): {\bigotimes_i}\, \lZ(X_1^i, A_1^i) \to {\bigotimes_i}\, \lZ(X_2^i, A_2^i)
	$$
\end{itemize}

\begin{rmk}
The apparently artificial definition of $\mathfrak{Z}$ can be better understood in terms of the corresponding map of $\lZ$. To see it, let $(W, A): (X_1, A_1) \to (X_2, A_2)$ and $(W', A'): (X_1', A_1') \to (X_2', A_2')$ be two tubes. Recall that, since $\lZ$ is lax monoidal, we have a natural trasformation $\Delta: \lZ(-) \otimes_R \lZ(-) \Rightarrow \lZ(- \sqcup -)$ (see section \ref{subsec:bimodules}). Then, $\mathfrak{Z}\left((W,A) \sqcup (W, A')\right)$ is just a lift
	\[
\begin{displaystyle}
   \xymatrix
   {
      \lZ(X_1, A_1) \otimes_R \lZ(X_1', A_1') \ar[d]_{\Delta_{(X_1, A_1), (X_1', A_1')}} \ar@{--{>}}[rrr]^{\mathfrak{Z}\left((W,A) \sqcup (W, A')\right)} &&& \lZ(X_2, A_2) \otimes_R \lZ(X_2', A_2') \ar[d]^{{\Delta_{(X_2, A_2), (X_2', A_2')}}} \\
   \lZ(X_1 \sqcup X_1', A_1 \cup A_1') \ar[rrr]_{\lZ\left((W,A) \sqcup (W, A')\right)} &&& \lZ(X_2 \sqcup X_2', A_2 \cup A_2')
   }
\end{displaystyle}   
\]
Since composition of tubes is performed by componentwise gluing, these lifts behave well with composition, which implies that $\mathfrak{Z}$ is well defined. This property no longer holds for $\lZ$ since there exists bordisms mixing several components, as the pair of pants for $n=2$.
\end{rmk}

In particular, from the TQFT described in Theorem \ref{thm:existence-s-LTQFT}, we obtain an almost-TQFT of pairs, $\mathfrak{Z}: \CTubpp{n} \to \Mod{R}$,
with $R = \K{\MHS{\QQ}}$, that allows the computation of Deligne-Hodge polynomials of representation varieties. To be precise, from the previous construction we obtain the following result.

\begin{cor}\label{cor:almost-TQFT}
There exists an almost-TQFT of pairs, $\mathfrak{Z}: \CTubpp{n} \to \Mod{\K{\MHS{\QQ}}}$, such that, for any $n$-dimensional connected closed orientable manifold $W$ and any finite set $A \subseteq W$ we have
$$
	\mathfrak{Z}(W,A)(\QQ_0) = \left[H_c^\bullet(\Rep{G}(W); \QQ)\right] \times \left[H_c^\bullet(G; \QQ)\right]^{|A|-1}
$$
where $\QQ_0 \in \K{\MHS{\QQ}}$ is the unit. In particular,
$$
	\DelHod{\mathfrak{Z}(W, A)(\QQ_0)} = \DelHod{G}^{|A|-1} \DelHod{\Rep{G}(W)}.
$$
\end{cor}

\begin{figure}[h]
	\begin{center}
	\includegraphics[scale=0.45]{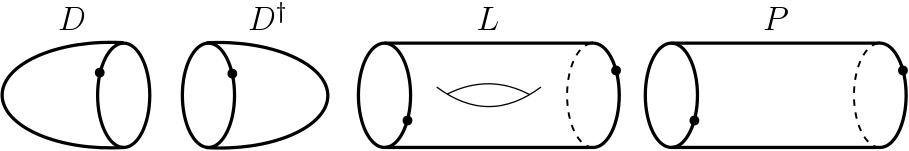}
	\caption{A set of generators of $\CTubpp{2}$.}
	\label{img:non-parabolic-tubes}
	\end{center}
\end{figure}

In the case $n=2$, we can go a step further and describe explicitly this functor. With respect to objects of $\CTubpp{2}$, recall that every non-empty $1$-dimensional closed manifold is diffeomorphic to $S^1$. Since $\pi_1(S^1)=\ZZ$ we have $\Rep{G}(S^1)=\Hom(\ZZ, G)=G$ and, thus, for any $x_0 \in S^1$
$$
	\mathfrak{Z}(\emptyset) = \K{\MHM{1}} = \K{\MHS{\QQ}}=R, \hspace{1cm} \mathfrak{Z}(S^1, \star) = \K{\MHM{G}}.
$$
Regarding morphisms, observe that, adapting the proof of Proposition 1.4.13 of \cite{Kock:2004} for the generators of $\CBordp{2}$, we have that a set of generators of $\CTubpp{2}$ is $\Delta = \left\{D, D^\dag, L, P\right\}$ where $D: \emptyset \to (S^1, \star)$ is the disc with a marked point in the boundary, $D^\dag: (S^1, \star) \to \emptyset$ is the opposite disc, $L: (S^1, \star) \to (S^1, \star)$ is the torus with two holes and a puncture on each boundary component, and $P: (S^1, \star) \to (S^1, \star)$ is the bordism $S^1 \times [0,1]$ with a puncture on each boundary component. They are depicted in Figure \ref{img:non-parabolic-tubes}.

For $D$ and $D^\dag$ the situation is simple since they are simply connected. Therefore, their images under the geometrisation functor $\Geo$ are
$$
\Geo(D) = \left[1 \longleftarrow 1 \stackrel{i}{\longrightarrow} G\right],
	 \hspace{1cm} \Geo(D^\dag) = \left[G \stackrel{i}{\longleftarrow} 1 \longrightarrow 1\right],
$$
so their images under $\mathfrak{Z}$ are
$$
	\mathfrak{Z}(D) = i_!: R=\K{\MHM{1}} \to \K{\MHM{G}} \hspace{1cm} \mathfrak{Z}(D^\dag) = i^*: \K{\MHM{G}} \to \K{\MHM{1}} = R.
$$

\begin{figure}[h]
	\begin{center}
	\includegraphics[scale=0.15]{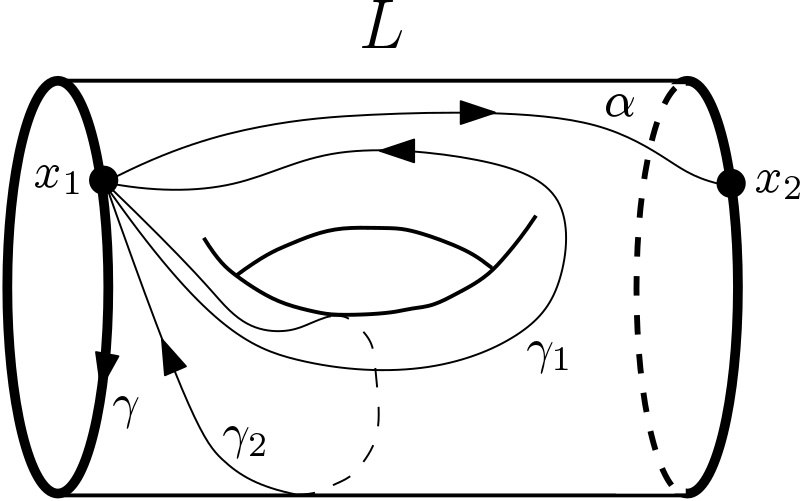}
	\caption{Chosen paths for $L$.}
	\label{img:paths-T}
	\end{center}
\end{figure}

For the holed torus $L: S^1 \to S^1$ the situation is a bit more complicated. Let $L=(T, A)$ where $A = \left\{x_1, x_2\right\}$ the set of marked points of $L$, $x_1$ in the ingoing boundary and $x_2$ in the outgoing boundary. Recall that $T$ is homotopically equivalent to a bouquet of three circles so its fundamental group is the free group with three generators. Thus, we can take $\gamma, \gamma_1, \gamma_2$ as the set of generators of $\Pi(L)_{x_1} = \pi_1(T, x_1)$ depicted in Figure \ref{img:paths-T} and $\alpha$ the shown path between $x_1$ and $x_2$.

With this description, $\gamma$ is a generator of $\pi_1(S^1, x_1)$ and $\alpha\gamma[\gamma_1, \gamma_2]\alpha^{-1}$ is a generator of $\pi_1(S^1, x_2)$, where $[\gamma_1, \gamma_2] = \gamma_1\gamma_2\gamma_1^{-1}\gamma_2^{-1}$ is the group commutator. Hence, since $\Rep{G}(L) =  \Hom(\Pi(T, A), G) = G^4$ we have that the geometrisation $\Geo(L)$ is the span
$$
\begin{matrix}
	G & \stackrel{p}{\longleftarrow} & G^4 & \stackrel{q}{\longrightarrow} & G \\
	g & \mapsfrom & (g, g_1, g_2, h) & \mapsto & hg[g_1,g_2]h^{-1}
\end{matrix}
$$
where $g, g_1, g_2$ and $h$ are the images of $\gamma, \gamma_1, \gamma_2$ and $\alpha$, respectively. Hence, we obtain that
$$
	\mathfrak{Z}(L): \K{\MHM{G}} \stackrel{p^*}{\longrightarrow} \K{\MHM{G^4}} \stackrel{q_!}{\longrightarrow} \K{\MHM{G}}.
$$

\begin{figure}[h]
	\begin{center}
	\includegraphics[scale=0.15]{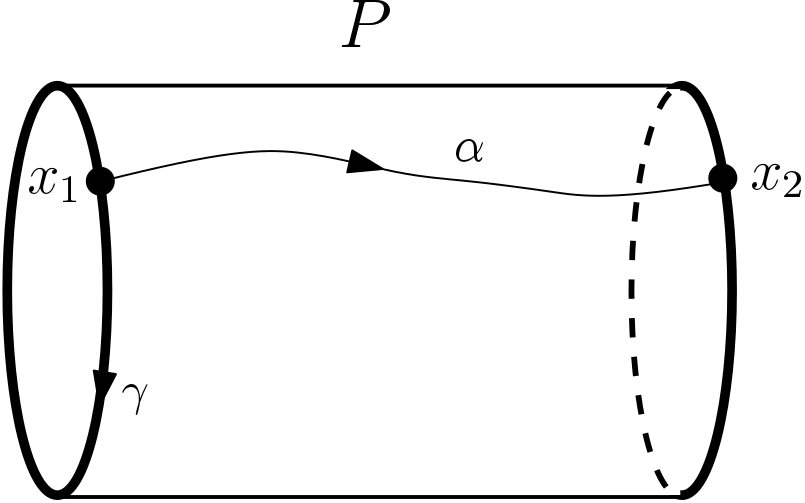}
	\caption{Chosen path for $P$.}
	\label{img:paths-P}
	\end{center}
\end{figure}

For the morphism $P$, let $P=(S^1 \times [0,1], A)$ where $A = \left\{x_1, x_2\right\}$ with $x_1, x_2$ the ingoing and outgoing boundary points respectively. Since $\pi_1(S^1 \times [0,1]) = \ZZ$, the fundamental groupoid $\Pi(P)$ has two vertices isomorphic to $\ZZ$ and, thus, $\Hom(\Pi(P),G)=G^2$. Let $\gamma$ be a generator of $\pi_1(S^1, x_1)$ and $\alpha$ a path between $x_1$ and $x_2$, as depicted in Figure \ref{img:paths-P}.

Since $\alpha_1\gamma\alpha_1^{-1}$ is a generator of $\pi_1(S^1, x_2)$ we obtain that $\Geo(P)$ is the span
$$
\begin{matrix}
	G & \stackrel{u}{\longleftarrow} & G^2 & \stackrel{v}{\longrightarrow} & G \\
	g & \mapsfrom & (g, h) & \mapsto & hgh^{-1}
\end{matrix}
$$
Hence, we have that
$$
	\mathfrak{Z}(P): \K{\MHM{G}} \stackrel{u^*}{\longrightarrow} \K{\MHM{G^3}} \stackrel{v_!}{\longrightarrow} \K{\MHM{G}}.
$$

Now, let $\Sigma_g$ be the closed oriented surface of genus $g$. If we choose any $g+1$ points we have a decomposition of the bordism $\Sigma_g: \emptyset \to \emptyset$ as $\Sigma_g = D^\dag \circ L^g \circ D$. Thus, by Corollary \ref{cor:almost-TQFT} we have that
$$
	\DelHod{\Rep{G}(\Sigma_g)} = \frac{1}{\DelHod{G}^g}\DelHod{\mathfrak{Z}(D^\dag) \circ \mathfrak{Z}(L)^g \circ \mathfrak{Z}(D) (\QQ_0)}.
$$
This kind of computations were carried out in the paper \cite{MM} for $G=\SL(2, \CC)$ and in \cite{Martinez:2017} for $G=\PGL(2,\CC)$. In future work, we shall explain the computations of those papers in these terms and extend them to the parabolic case using the techniques explained in the following section.

\subsection{The parabolic case}

The previous constructions can be easily extended to the parabolic case. Following the construction above, it is just necessary to adapt the geometrisation functor to the parabolic context.

As above, let us fix a complex algebraic group $G$ and, as parabolic data, we choose for $\Lambda$ a collection of subsets of $G$ closed under conjugacy. Then, as geometrisation functor, we define the $2$-functor $\Geo: \CBordpp{n}(\Lambda) \to \cS = \Span{\CVar}$ as follows:
\begin{itemize}
	\item For any object $(X, A)$ of $\CBordpp{n}(\Lambda)$, we set $\Geo(X,A)= \Rep{G}(X, A) = \Hom(\Pi(X, A), G)$.
	\item Let $(W, A, Q): (X_1, A_1) \to (X_2, A_2)$ with $Q = \left\{(S_1, \lambda_1), \ldots, (S_s, \lambda_s)\right\}$ be a $1$-morphism. For short, let us denote $S = S_1 \cup \ldots \cup S_s$. To this morphism we assign the span
	$$
		\Rep{G}(X_2, A_2) \longleftarrow \Rep{G}(W, A, Q) \longrightarrow \Rep{G}(X_1, A_1)
	$$
where $\Rep{G}(W, A, Q) \subseteq \Rep{G}(W - S, A)$ is the subvariety of groupoid homomorphisms $\rho: \Pi(W - S, A) \to G$ satisfying the following condition: if $\gamma$ is a loop of $\pi_1(W - S, a)$, being $a$ any point of $W$ in the connected component of $\gamma$, whose image under the map induced by inclusion $\pi_1(W - S, a) \to \pi_1((W - S) \cup S_k, a)$ vanishes, then $\rho(\gamma) \in \lambda_k$. Observe that, since the $\lambda \in \Lambda$ are closed under conjugation, that condition does not depend on the chosen basepoint $a$.
	\item For a $2$-morphism $(W, A, Q) \Rightarrow (W', A', Q')$ given by a diffeomorphism $F: W \to W'$ we use the induced map $F^*: \Rep{G}(W',A',Q') \to \Rep{G}(W,A,Q)$ to create a commutative diagram
	\[
\begin{displaystyle}
   \xymatrix
   {	& \Rep{G}(W', A', Q') \ar[rd] \ar[ld] \ar[d]^{F^*} & \\
   		\Rep{G}(X_2, A_2) & \Rep{G}(W, A, Q)  \ar[r] \ar[l] & \Rep{G}(X_1, A_1)
   }
\end{displaystyle}   
\]
which is a $2$-cell in $\Span{\CVar}$.
\end{itemize}
As in section \ref{subsec:E-pol-geometric-part}, these assigments can be put together to define a $2$-functor $\Geo: \CBordpp{n}(\Lambda) \to \cS$. This functor, together with the algebraisations $\AlgA: \CVar \to \Rng$ and $\AlgB: \CVar \to \Mod{R}$ described in section \ref{subsec:E-pol-algebraic-part}, gives us a parabolic soft TQFT, $\sZ: \CBordpp{n}(\Lambda) \to \Bim{R}$ and a parabolic lax monoidal TQFT, $\lZ: \CBordpp{n}(\Lambda) \to \Mod{R}_t$. As in the non-parabolic case, the functor $\lZ$ satisfies that, for any closed connected orientable $n$-dimensional manifold $X$ any finite set $A \subseteq X$ and any parabolic structure $Q$ on $X$ disjoint from $A$ we have
$$
	\lZ(X, A, Q)(\QQ_0) = \left[H_c^\bullet\left(\Rep{G}(X, Q); \QQ\right)\right] \otimes \left[H_c^\bullet\left(G; \QQ\right)\right]^{|A|-1}.
$$
As always, this implies that
$$
	\DelHod{\lZ(X, A, Q)(\QQ_0)} = \DelHod{G}^{|A|-1} \DelHod{\Rep{G}(X, Q)}.
$$ 

Furthermore, following the construction in section \ref{subsec:almost-TQFT}, we can also modify this result to obtain an almost-TQFT, $\mathfrak{Z}: \CTubpp{n}(\Lambda) \to \Mod{\K{\MHS{\QQ}}}$, where $\CTubpp{n}(\Lambda)$ is the category of tubes of pairs with parabolic data $\Lambda$. In this parabolic case, for any closed orientable $n$-dimensional manifold $W$ with parabolic structure $Q$ and any finite set $A \subseteq W$ we obtain
\begin{align*}
	\mathfrak{Z}(W, A, Q)(\QQ_0) &= \left[H_c^\bullet(\Rep{G}(W, Q); \QQ)\right] \times \left[H_c^\bullet(G; \QQ)\right]^{|A|-1}, \\
\DelHod{\mathfrak{Z}(W, A, Q)(\QQ_0)} &= \DelHod{G}^{|A|-1} \DelHod{\Rep{G}(W, Q)}.
\end{align*}

\begin{figure}[h]
	\begin{center}
	\includegraphics[scale=0.15]{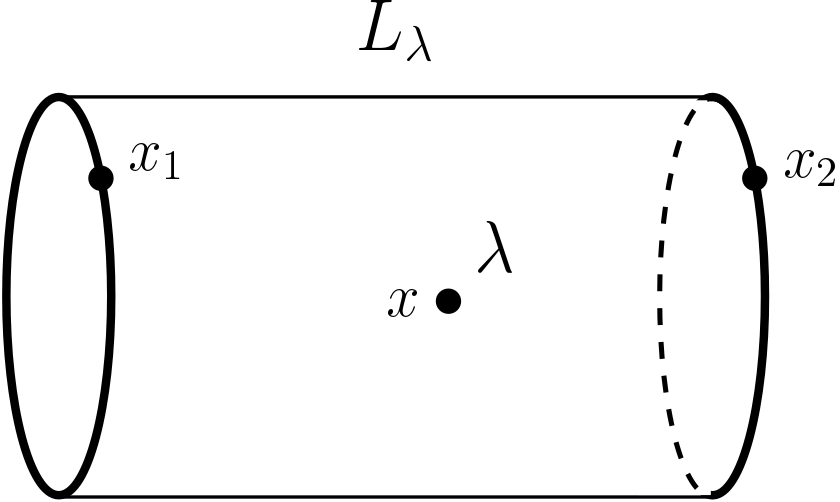}
	\caption{Tube with marked point.}
	\label{img:parabolic-tube}
	\end{center}
\end{figure}

In the particular case $n=2$, observe that $\Par{W}$ is just the set of interior points of the surface $W$. Therefore, in order to obtain a set of generators of $\CTubpp{n}(\Lambda)$ it is enough to consider the elements of the set of generators $\Delta$ previouly defined with no parabolic structure and to add a collection of tubes $L_\lambda = (S^1 \times [0,1], \left\{x_1, x_2\right\}, \left\{(x, \lambda)\right\})$ for $\lambda \in \Lambda$, where $x$ is any interior point of $S^1 \times [0,1]$ with parabolic structure $\lambda$ and $x_1, x_2$ are points in the ingoing and outgoing boundaries repectively, as depicted in Figure \ref{img:parabolic-tube}.

In this case, observe that $\pi_1((S^1 \times [0,1]) - \left\{x_0\right\})$ is the free group with two generators and that the fundamental groupoid of $L_\lambda$ has two vertices. These generators can be taken to be around the incoming boundary and around the marked point so $\Rep{G}(L_\lambda) = G^2 \times \lambda$. Thus, the geometrisation $\Geo(L_\lambda)$ is the span
$$
\begin{matrix}
	G & \stackrel{r}{\longleftarrow} & G^2 \times \lambda & \stackrel{s}{\longrightarrow} & G \\
	g & \mapsfrom & (g, g_1, h) & \mapsto & g_1ghg_1^{-1}
\end{matrix}
$$
Hence, the image $\mathfrak{Z}(L_\lambda)$ is the morphism
$$
	\mathfrak{Z}(L_\lambda): \K{\MHM{G}} \stackrel{r^*}{\longrightarrow} \K{\MHM{G^2 \times \lambda}} \stackrel{s_!}{\longrightarrow} \K{\MHM{G}}.
$$
With this description we have proven the following result.

\begin{thm}\label{thm:almost-tqft-parabolic}
Let  $\Sigma_g$ be a closed oriented surface of genus $g$ and $Q$ a parabolic structure on $\Sigma_g$ with $s$ marked points with data $\lambda_1, \ldots, \lambda_s \in \Lambda$,
then, 
$$
	\DelHod{\Rep{G}(\Sigma_g, Q)} = \frac{1}{\DelHod{G}^{g+s}}\DelHod{\mathfrak{Z}(D^\dag) \circ \mathfrak{Z}(L_{\lambda_s}) \circ \ldots \circ \mathfrak{Z}(L_{\lambda_1}) \circ \mathfrak{Z}(L)^g \circ \mathfrak{Z}(D) (\QQ_0)}.
$$
\end{thm}

This formula gives a general recipe for computing Deligne-Hodge polynomials of representation varieties, for any group $G$. That is, once computed explicitly the homomorphism $\mathfrak{Z}(L): \K{\MHM{G}} \to \K{\MHM{G}}$, all the Deligne-Hodge polynomials are known in the non-parabolic case and, if we also compute $\mathfrak{Z}(L_\lambda): \K{\MHM{G}} \to \K{\MHM{G}}$ for any conjugacy class $\lambda \subseteq G$, also the parabolic case follows.

\section{Some examples of computations}

In this section, we will discuss some toy examples of the previous TQFT in order to show how it works and how it can be used to give an effective method of computation.

\subsection{Finite groups} \label{subsec:finite-group-calculation}

Suppose that $X$ is a complex algebraic variety that can be decomposed as $X = Z \sqcup U$, with $i: Z \hookrightarrow X$ Zariski closed and $j: U \hookrightarrow X$ open. Then, in Section 4.4 of \cite{Saito:1990}, it is proven that $i_! + j_!: \K{\MHM{Z}} \otimes \K{\MHM{U}} \to \K{\MHM{X}}$ and $i^* \oplus j^*: \K{\MHM{X}} \to \K{\MHM{Z}} \otimes \K{\MHM{U}}$ are isomorphisms, and $i_!i^* + j_!j^*: \K{\MHM{X}} \to \K{\MHM{X}}$ is the identity morphism. In particular, if $X$ is finite, we have an isomorphism
$$
	\sum_{x \in X} (i_x)_!: \bigoplus_{x \in X} \K{\MHM{x}}=\K{\MHS{\QQ}}^{|X|} \to \K{\MHM{X}},
$$
where $i_x: \set{x} \hookrightarrow X$ is the inclusion of $x \in X$ and $|X|$ denotes the cardinal of $X$. This implies that $\K{\MHM{X}}$ is the free $\K{\MHS{\QQ}}^{|X|}$-module generated by the $|X|$ generators $\QQ_x = (i_x)_! {\QQ}_0$. Moreover, if $f: Y \to X$ is a regular morphism between finite varieties, then
$$
	f_! \underline{\QQ}_Y = \sum_{x \in X} |f^{-1}(x)|\,\QQ_x.
$$

As an application, suppose that $G$ is a finite group with $n$ elements. In that case, we have an identification $\K{\MHM{G}} = \K{\MHS{\QQ}}^n$. Moreover, the image under the TQFT of the disc is $\mathfrak{Z}(D): \K{\MHS{\QQ}} \to \K{\MHS{\QQ}}^n$ given by $\mathfrak{Z}(D)(\QQ_0) = \QQ_1$ with $1 \in G$ the unit. On the other hand, $\mathfrak{Z}(D^\dag): \K{\MHS{\QQ}}^n \to \K{\MHS{\QQ}}$ is nothing but the projection onto $\QQ_1$.

With respect to the holed torus $L: (S^1, \star) \to (S^1, \star)$, fix $g \in G$. Then, we have a commutative diagram
\[
\begin{displaystyle}
   \xymatrix
   {	& G^3 \ar[r] \ar[d] \ar[ld]_{\pi_g} & g \ar[d]^{i_g} \\
   		G & G^4 \ar[r]_{p} \ar[l]^q & G
      }
\end{displaystyle}   
\]
where the leftmost arrow is given by $(g_1, g_2, h) \mapsto (g, g_1, g_2,h)$ and $\pi_g(g_1, g_2, h) = hg[g_1,g_2]h^{-1}$. Observe that the rightmost square of this diagram is a pullback. Hence, by the Beck-Chevalley condition and the previous formula we have that
\begin{align*}
\mathfrak{Z}(L)(\QQ_g) &= q_!p^*(i_g)_!\QQ_0 = (\pi_g)_!\underline{\QQ}_{G^3} = \sum_{a \in G} \left|\left\{(g_1, g_2, h) \in G^3\,|\, hg[g_1,g_2]h^{-1}=a\right\}\right|\,\QQ_a.
\end{align*}
Using the previous description of $\K{\MHM{G}}$, this formula fully determines $\mathfrak{Z}(L)$. Analogous considerations can also be done for the parabolic tubes $L_\lambda: (S^1, \star) \to (S^1, \star)$. In this way, in the finite setting, the almost-TQFT is just a systematic way of counting the number of points in the representation variety.

\subsection{General affine group} \label{subsec:general-affine-calculation}

Let us take $G = \Aff{\CC}$, the group of $\CC$-linear affine transformations of the complex line. Recall that its elements are the matrices of the form
$
\begin{pmatrix}
	a & b\\
	0 & 1\\
\end{pmatrix}
$,
with $a \in \CC^* = \CC - \set{0}$ and $b \in \CC$. The group operation is given by matrix multiplication. In this way, $\Aff{\CC}$ is isomorphic to the semidirect product $\CC^* \ltimes_\varphi \CC$ with the action $\varphi: \CC^* \times \CC \to \CC$, $\varphi(a,b)=ab$.

In order to compute the morphism $\mathfrak{Z}(L): \K{\MHM{\Aff{\CC}}} \to \K{\MHM{\Aff{\CC}}}$, recall that, with the notation of Section \ref{subsec:almost-TQFT}, $\mathfrak{Z}(L) = q_!p^*$ and $\mathfrak{Z}(D)=i_!$. We have a commutative diagram
	\[
\begin{displaystyle}
   \xymatrix
   {	& \Aff{\CC}^3 \ar[r]^c \ar[d] \ar[ld]_{\varpi} & 1 \ar[d]^{i} \\
   		\Aff{\CC} & \Aff{\CC}^4 \ar[r]_{p} \ar[l]^q & \Aff{\CC}
      }
\end{displaystyle}   
\]
where the leftmost vertical arrow is given by $(A_1, A_2, B) \mapsto (I, A_1, A_2, B)$ and $\varpi(A_1, A_2, B) = B[A_1, A_2]B^{-1}$, being $I \in \Aff{\CC}$ the unit matrix. Moreover, the square is a pullback, so we have $\mathfrak{Z}(L) \circ \mathfrak{Z}(D)(\QQ_1) = q_!p^*i_!\QQ_0 = \varpi_! c^*\QQ_0 = \varpi_! \underline{\QQ}_{\Aff{\CC}^3}$. Explicitly, the morphism $\varpi$ is given by
$$
	\varpi\left(\begin{pmatrix}
	a_1 & b_1\\
	0 & 1\\
\end{pmatrix}, \begin{pmatrix}
	a_2 & b_2\\
	0 & 1\\
\end{pmatrix}, \begin{pmatrix}
	x & y\\
	0 & 1\\
\end{pmatrix}\right) = \begin{pmatrix}
	1 & (a_1-1)b_2x - (a_2-1)b_1x\\
	0 & 1\\
\end{pmatrix}.
$$
Therefore, $\varpi$ is a projection onto $\ASO{\CC} \subseteq \Aff{\CC}$, the subgroup of affine transformations whose linear part is orthogonal and orientation preserving. Outsite $I \in \ASO{\CC}$, $\varpi$ is a locally trivial fibration in the Zariski topology with fiber, on $\alpha \neq 0$, given by
\begin{align*}
F &= \left\{(a_1, a_2, x, b_1, b_2, y) \in (\CC^*)^3 \times \CC^3\,|\,(a_1-1)b_2x - (a_2-1)b_1x = \alpha\right\} \\
	& = \left\{b_2 = \frac{\alpha + (a_2-1)b_1x}{(a_1-1)x}, a_1 \neq 1\right\} \sqcup \left\{b_1 = -\frac{\alpha}{(a_2-1)x}, a_1 = 1\right\} \\
	& \cong \left[\left(\CC-\left\{0,1\right\}\right) \times (\CC^*)^2 \times \CC^2 \right] \sqcup \left[\CC-\set{0,1} \times \CC^* \times \CC^2 \right].
\end{align*}
On the other hand, on the identity matrix $I$, the special fiber is
\begin{align*}
	\varpi^{-1}(I) &= \left\{(a_1, a_2, x, b_1, b_2, y) \in (\CC^*)^3 \times \CC^3\,|\,(a_1-1)b_2 = (a_2-1)b_1\right\} \\
	& = \left\{b_2 = \frac{(a_2-1)b_1}{a_1-1}, a_1 \neq 1\right\} \sqcup \left\{a_1 = 1, a_2=1\right\} \sqcup \left\{a_1 = 1, a_2 \neq 1, b_1 = 0\right\} \\
	& \cong \left[\left(\CC-\left\{0,1\right\}\right) \times (\CC^*)^2 \times \CC^2 \right] \sqcup \left[\CC^* \times \CC^3\right] \sqcup \left[\CC-\set{0,1} \times \CC^* \times \CC^2 \right].
\end{align*}

In \cite{Saito:1990} (in particular, Theorem 3.27 and Remark on p.313), it is shown that if $\pi: X \to B$ is a regular morphism of complex varieties with trivial monodromy and fiber $F$, then $\pi_! \underline{\QQ}_X = [H_c^\bullet(F; \QQ)]\underline{\QQ}_B$. Here, $[H_c^\bullet(F; \QQ)] \in \K{\MHS{\QQ}}$ is the image, in the Grothendieck ring, of the mixed Hodge structure on the cohomology (with compact support) of $F$.

As an application, let us denote $\ASO{\CC}^* = \ASO{\CC} - \set{I}$ with inclusion $j: \ASO{\CC}^* \hookrightarrow \Aff{\CC}$. Then, we have that
$$
	\varpi_!\underline{\QQ} = i_!i^*\varpi_! \underline{\QQ} + j_!j^*\varpi_!\underline{\QQ} =  i_!(\varpi|_{\varpi^{-1}(I)})_!\underline{\QQ} + j_!(\varpi|_{\varpi^{-1}(\ASO{\CC}^*)})_! \underline{\QQ}.
$$

For the first map, recall that $\varpi$ is locally trivial in the Zariski topology over $\ASO{\CC}$. Using the previous computation, we also have that the Hodge structure on the fiber is $[H_c^\bullet(F)] = (q-2)(q-1)^2q^2 + (q-2)(q-1)q^2 = q(q-1)(q^3-2q^2)$, where $q = [\QQ(-1)] = [H_c^\bullet(\CC)]$ is the Tate Hodge structure of weight $2$. On the other hand, the Hodge structure of the special fiber is $[H_c^\bullet(\varpi^{-1}(I))] = (q-2)(q-1)^2q^2 + (q-1)q^3 + (q-2)(q-1)q^2 = q(q-1)(q^3-q^2)$.
Hence, putting all together, we obtain that
\begin{align*}
	\mathfrak{Z}(L) \circ \mathfrak{Z}(D)(\QQ_0) &=  i_!(\varpi|_{\varpi^{-1}(I)})_!\underline{\QQ} + j_!(\varpi|_{\varpi^{-1}(\ASO{\CC}^*)})_! \underline{\QQ} \\
	&= q(q-1)(q^3-q^2) \,i_!\QQ_0 + q(q-1)(q^3-2q^2)\,j_!\underline{\QQ}_{\ASO{\CC}^*}.
\end{align*}

In this way, if we want to apply $\mathfrak{Z}(L)$ twice, we need to compute the image $\mathfrak{Z}(L)(j_!\underline{\QQ}_{\ASO{\CC}^*})$. This computation is quite similar to the previous one. Firstly, we again have a commutative diagram whose square is a pullback
	\[
\begin{displaystyle}
   \xymatrix
   {	& \ASO{\CC}^* \times \Aff{\CC}^3 \ar[r] \ar[d] \ar[ld]_{\vartheta} & \ASO{\CC}^* \ar[d]^{j} \\
   		\Aff{\CC} & \Aff{\CC}^4 \ar[r]_{p} \ar[l]^q & \Aff{\CC}
      }
\end{displaystyle}  
\]
Here, the leftmost vertical arrow is the inclusion map and $\vartheta(A, A_1, A_2, B) = BA[A_1, A_2]B^{-1}$. Computing explicitly, we have that
$$
	\vartheta\left(\begin{pmatrix}
	1 & \beta\\
	0 & 1\\
\end{pmatrix}, \begin{pmatrix}
	a_1 & b_1\\
	0 & 1\\
\end{pmatrix}, \begin{pmatrix}
	a_2 & b_2\\
	0 & 1\\
\end{pmatrix}, \begin{pmatrix}
	x & y\\
	0 & 1\\
\end{pmatrix}\right) = \begin{pmatrix}
	1 & (a_1-1)b_2x - (a_2-1)b_1x + \beta x\\
	0 & 1\\
\end{pmatrix}.
$$
Hence, $\vartheta$ is again a morphism onto $\ASO{\CC} \subseteq \Aff{\CC}$. Over $I \in \ASO{\CC}$, the fiber is
\begin{align*}
	\vartheta^{-1}(I) &= \left\{(\beta, a_1, a_2, x, b_1, b_2, y) \in (\CC^*)^4 \times \CC^3\,|\, (a_2-1)b_1x - (a_1-1)b_2x = \beta\right\} \\
	&= \left[(\CC^*)^3 \times \CC^3\right]- \left\{(a_1-1)b_2 - (a_2-1)b_1 = 0\right\} = \left[(\CC^*)^3 \times \CC^3\right]-\varpi^{-1}(I).
\end{align*}
Thus, $[H_c^\bullet(\vartheta^{-1}(I))] = (q-1)^3q^3-q(q-1)(q^3-q^2) = q(q-1)(q^4-3q^3+2q^2)$.

Analogously, on $\ASO{\CC}^*$, we have that $\vartheta$ is a locally trivial fibration in the Zariski topology with fiber over $\alpha \neq 0$
\begin{align*}
	F' &= \left\{(\beta, a_1, a_2, x, b_1, b_2, y) \in (\CC^*)^4 \times \CC^3\,|\,(a_1-1)b_2x - (a_2-1)b_1x + \beta = \alpha\right\} \\
	&= \left[(\CC^*)^3 \times \CC^3\right]- \left\{(a_1-1)b_2 - (a_2-1)b_1 = \alpha\right\} = \left[(\CC^*)^3 \times \CC^3\right]-F.
\end{align*}
Hence, the Hodge structure on the cohomology of the fiber is $[H_c^\bullet(F')] = (q-1)^3q^3-q(q-1)(q^3-2q^2) = q(q-1)(q^4-3q^3+3q^2)$. Putting together these computations we obtain that
\begin{align*}
	\mathfrak{Z}(L)\left(j_!\underline{\QQ}_{\ASO{\CC}^*}\right) &= \vartheta_!\underline{\QQ} =  i_!(\vartheta|_{\vartheta^{-1}(I)})_!\underline{\QQ} + j_!(\vartheta|_{\vartheta^{-1}(\ASO{\CC}^*)})_! \underline{\QQ} \\
	&= q(q-1)(q^4-3q^3+2q^2) \,i_!\QQ_0 + q(q-1)(q^4-3q^3+3q^2)\,j_!\underline{\QQ}_{\ASO{\CC}^*}.
\end{align*}

Let $W \subseteq \K{\MHM{\Aff{\CC}}}$ be the submodule generated by the elements $i_!\QQ_0$ and $j_!\underline{\QQ}_{\ASO{\CC}^*}$. The previous computation has shown that $\mathfrak{Z}(L)(W) \subseteq W$. Furthermore, indeed we have $W = \langle \mathfrak{Z}(L)^g(i_!{\QQ}_0)\rangle_{g=0}^\infty$. On $W$, the morphism $\mathfrak{Z}(D^\dag): W \to \K{\MHS{\QQ}}$ is given by the projection $\mathfrak{Z}(D^\dag)(i_!\QQ_0) = \QQ_0$ and $\mathfrak{Z}(D^\dag)(j_!\underline{\QQ}_{\ASO{\CC}^*}) = 0$. Hence, regarding the computation of Deligne-Hodge polynomials of representation varieties, we can restrict our attention to $W$.

If we want to explicitly compute these polynomials, observe that, by the previous calculations, in the set of generators $i_!\QQ_0, j_!\underline{\QQ}_{\ASO{\CC}^*}$ of $W$, the matrix of $\mathfrak{Z}(L): W \to W$ is
$$
\mathfrak{Z}(L) = q(q-1) \begin{pmatrix}
	q^3-q^2 & q^4-3q^3+2q^2 \\
	q^3-2q^2 & q^4-3q^3+3q^2\\
\end{pmatrix}.
$$

Let us denote $q = uv \in \ZZ[u^{\pm 1}, v^{\pm}]$. This abuse of notation is compatible with the fact that $\DelHod{\CC} = \DelHod{q}=uv$. Since $\DelHod{\Aff{\CC}} = \DelHod{\CC^* \times \CC} = q(q-1)$, using the formula of Theorem \ref{thm:almost-tqft-parabolic} we obtain that
\begin{align*}
	\DelHod{\Rep{\Aff{\CC}}(\Sigma_g)} &= \begin{pmatrix}
	1 & 0\\
\end{pmatrix}
\begin{pmatrix}
	q^3-q^2 & q^4-3q^3+2q^2 \\
	q^3-2q^2 & q^4-3q^3+3q^2\\
\end{pmatrix}^g\begin{pmatrix}
	1\\
	0\\
\end{pmatrix} \\
&= \begin{pmatrix}
	1 & 0\\
\end{pmatrix}
\begin{pmatrix}
	q-1 & q-1 \\
	-1 & q-1\\
\end{pmatrix}
\begin{pmatrix}
	q^{2g} & 0 \\
	 & q^{2g}(q-1)^{2g}\\
\end{pmatrix}\begin{pmatrix}
	q-1 & q-1 \\
	-1 & q-1\\
\end{pmatrix}^{-1}\begin{pmatrix}
	1\\
	0\\
\end{pmatrix} \\
&=q^{2g-1}\left((q-1)^{2g}+q-1\right).
\end{align*}

\begin{rmk}
This formula can be checked directly since, in this toy model, a recursive formula for the Deligne-Hodge polynomial can be easily found. Observe that we have an explicit expression of the product of commutators
$$
	\prod_{i=1}^g \left[\begin{pmatrix}
	a_{2i-1} & b_{2i-1} \\
	0 & 1\\
\end{pmatrix},\begin{pmatrix}
	a_{2i} & b_{2i} \\
	0 & 1\\
\end{pmatrix}\right] = \begin{pmatrix}
	1 & {\displaystyle \sum_{i=1}^g (a_{2i-1}-1)b_{2i} - (a_{2i}-1)b_{2i-1} }\\
	0 & 1\\
\end{pmatrix}.
$$
Therefore, if we define the auxiliar variety
$$
	X_{k} = \left\{(a_1, b_1, \ldots, a_{k}, b_{k}) \in (\CC^* \times \CC)^k\,\left|\, \sum_{i=1}^{k} (a_i-1)b_i = 0\right.\right\},
$$
we have that $\Rep{\Aff{\CC}}(\Sigma_g) \cong X_{2g}$. The important point is that the varieties $X_k$ can be recursively computed as
\begin{align*}
	X_k &= \left\{b_k = -\frac{\sum_{i=1}^{k-1} (a_i-1)b_i}{a_k-1}, a_k \neq 1\right\} \sqcup \left\{\sum_{i=1}^{k-1} (a_i-1)b_i,a_k=1\right\} \\
	&= \left[(\CC-\set{0,1}) \times \left(\CC^* \times \CC\right)^{k-1}\right] \sqcup \left[\CC \times X_{k-1}\right],
\end{align*}
together with the base case $X_1 = (\CC-\set{0,1}) \sqcup \CC$. This gives rise to the recursive formula for the Deligne-Hodge polynomials $\DelHod{X_k} = (q-2)q^{k-1}(q-1)^{k-1} + q\DelHod{X_{k-1}}$ with base case $\DelHod{X_1} = 2q-2$. Therefore, we obtain the formula $\DelHod{\Rep{\Aff{\CC}}(\Sigma_g)} = q^{2g}(q-1)^{2g-2}(q-2) + q^2\DelHod{\Rep{\Aff{\CC}}(\Sigma_{g-1})}$ with $\DelHod{\Rep{\Aff{\CC}}(\Sigma_1)} = q^3-q^2$, that agrees with the previous computation.
\end{rmk}

\subsection{Concluding remarks}

The previous examples strongly suggest that, despite that $\K{\MHM{G}}$ might be infinitely generated, all the relevant data for the computations are contained in a finitely generated module, at least in the non-parabolic case.

To be precise, given a complex algebraic group $G$, let $W_G = \langle \mathfrak{Z}(L)^g(i_!{\QQ}_0)\rangle_{g=0}^\infty \subseteq \K{\MHM{G}}$, the submodule generated by the image of the unit under all the closed surfaces with connected boundary. By definition, we have that $\mathfrak{Z}(L)(W_G) \subseteq W_G$ so, for computational purposes, we can restrict our attention to $W_G$.

It may be expected $W_G$ to be significatly smaller than the whole $\K{\MHM{G}}$. For example, in the example of Section \ref{subsec:general-affine-calculation}, we have that $W_G$ is generated by only two elements. Moreover, computations of \cite{MM}, in the case $G=\SL(2,\CC)$, and \cite{Martinez:2017} for $G=\PGL(2,\CC)$, show that $W_{G}$ is, in both cases, a finitely generated module (of $8$ and $6$ generators, respectively).
That is important because, in that case, the knowledge of the image of finitely many elements of $W_G$ is enough to characterize completely the map along $W_G$, as it happens in Sections \ref{subsec:finite-group-calculation} and \ref{subsec:general-affine-calculation}.
Therefore, there are good reasons to expect the following result, at least for some general class of algebraic groups in the non-parabolic case.

\begin{conjecture}
The module $W_G \subseteq \K{\MHM{G}}$ is finitely dimensional.
\end{conjecture}

Another related work consists of computing these polynomials explicitly for some particular groups $G$. In this direction, the upcoming papers \cite{GP-2018a} and \cite{GP-2018b} (see also \cite{Gonzalez-Prieto:Thesis}) perform this task for $G = \SL(2,\CC)$ and parabolic structures, a case that still remains unknown. In \cite{GP-2018a}, the almost-TQFT is computed in this case following the lines of Section \ref{subsec:general-affine-calculation}. However, in contrast with that section, the obtained fibrations have non-trivial monodromy that must be traced, making the calculation much more involved. Using this method, in the end, explicit expressions for the Deligne-Hodge polynomials of the parabolic $\SL(2,\CC)$-representation varieties are obtained. 

However, if we want to compute the corresponding polynomials for character varieties, we also need to analyze the GIT quotient under the conjugacy action. This is the gap filled by \cite{GP-2018b}, in which new stratification techniques for the quotient are developed. As a result, the Deligne-Hodge polynomials of these parabolic character varieties are computed.

Finally, the framework developed in this paper might be useful for understanding the mirror-symmetry conjectures related to character varieties, as stated in \cite{Hausel:2005} and \cite{Hausel-Rodriguez-Villegas:2008}. For example, computations of \cite{MM} and \cite{Martinez:2017} for the Langlands dual groups $\SL(2,\CC)$ and $\PGL(2,\CC)$ show that $W_{\SL(2,\CC)}$ and $W_{\PGL(2,\CC)}$ are strongly related, in the sense that their generators agree in the space of traces. Based on these observations, it would be interesting to study whether the linear maps $\mathfrak{Z}(L)$ are somehow related for Langland dual groups.

%%%%%%%%%%%%%%%%%%%%%%%%%%%%%%%%%%%%%%%%%%%%%%%%%%%%%%%%%%%%%%%%
%%%%%%%%%%%%%%%%%%%%%%% BIBLIOGRAPHY %%%%%%%%%%%%%%%%%%%%%%%%%%%
%%%%%%%%%%%%%%%%%%%%%%%%%%%%%%%%%%%%%%%%%%%%%%%%%%%%%%%%%%%%%%%%

%\nocite{*}
\bibliography{bibliography.bib}{}

\begin{thebibliography}{10}

\bibitem{Atiyah:1988}
M.~Atiyah.
\newblock Topological quantum field theories.
\newblock {\em Inst. Hautes \'Etudes Sci. Publ. Math.}, (68):175--186 (1989),
  1988.

\bibitem{Ayoub}
J.~Ayoub.
\newblock Les six op\'erations de {G}rothendieck et le formalisme des cycles
  \'evanescents dans le monde motivique. {I}.
\newblock {\em Ast\'erisque}, (314):x+466 pp. (2008), 2007.

\bibitem{Baez-Dolan}
J.~C. Baez and J.~Dolan.
\newblock Higher-dimensional algebra and topological quantum field theory.
\newblock {\em J. Math. Phys.}, 36(11):6073--6105, 1995.

\bibitem{Baraglia-Hekmati:2016}
D.~Baraglia and P.~Hekmati.
\newblock Arithmetic of singular character varieties and their
  {$E$}-polynomials.
\newblock {\em Proc. Lond. Math. Soc. (3)}, 114(2):293--332, 2017.

\bibitem{Beilinson-Drinfeld}
A.~A. Beilinson and V.~G. Drinfeld.
\newblock Quantization of {H}itchin's fibration and {L}anglands' program.
\newblock In {\em Algebraic and geometric methods in mathematical physics
  ({K}aciveli, 1993)}, volume~19 of {\em Math. Phys. Stud.}, pages 3--7. Kluwer
  Acad. Publ., Dordrecht, 1996.

\bibitem{Ben-Zvi-Francis-Nadler:2010}
D.~Ben-Zvi, J.~Francis, and D.~Nadler.
\newblock Integral transforms and {D}rinfeld centers in derived algebraic
  geometry.
\newblock {\em J. Amer. Math. Soc.}, 23(4):909--966, 2010.

\bibitem{Ben-Zvi-Gunningham-Nadler}
D.~Ben-Zvi, S.~Gunningham, and D.~Nadler.
\newblock The character field theory and homology of character varieties.
\newblock {\em Preprint arXiv:1705.04266}, 2017.

\bibitem{Ben-Zvi-Nadler}
D.~Ben-Zvi and D.~Nadler.
\newblock The character theory of a complex groups.
\newblock {\em Preprint arXiv:0904.1247}, 2009.

\bibitem{Ben-Zvi-Nadler:2016}
D.~Ben-Zvi and D.~Nadler.
\newblock Betti geometric langlands.
\newblock {\em Preprint arXiv:1606.08523}, 2016.

\bibitem{Benabou}
J.~B\'enabou.
\newblock Introduction to bicategories.
\newblock In {\em Reports of the {M}idwest {C}ategory {S}eminar}, pages 1--77.
  Springer, Berlin, 1967.

\bibitem{Biquard-Jardim}
O.~Biquard and M.~Jardim.
\newblock Asymptotic behaviour and the moduli space of doubly-periodic
  instantons.
\newblock {\em J. Eur. Math. Soc. (JEMS)}, 3(4):335--375, 2001.

\bibitem{Boden-Yokogawa}
H.~U. Boden and K.~Yokogawa.
\newblock Moduli spaces of parabolic {H}iggs bundles and parabolic {$K(D)$}
  pairs over smooth curves. {I}.
\newblock {\em Internat. J. Math.}, 7(5):573--598, 1996.

\bibitem{Brown:1967}
R.~Brown.
\newblock Groupoids and van {K}ampen's theorem.
\newblock {\em Proc. London Math. Soc. (3)}, 17:385--401, 1967.

\bibitem{Brown}
R.~Brown.
\newblock {\em Topology and groupoids}.
\newblock BookSurge, LLC, Charleston, SC, 2006.
\newblock Third edition of {{\i}t Elements of modern topology} [McGraw-Hill,
  New York, 1968; MR0227979], With 1 CD-ROM (Windows, Macintosh and UNIX).

\bibitem{Carlsson-Rodriguez-Villegas}
E.~Carlsson and F.~Rodriguez~Villegas.
\newblock Vertex operators and character varieties.
\newblock {\em Adv. Math.}, 330:38--60, 2018.

\bibitem{Corlette:1988}
K.~Corlette.
\newblock Flat {$G$}-bundles with canonical metrics.
\newblock {\em J. Differential Geom.}, 28(3):361--382, 1988.

\bibitem{DeligneII:1971}
P.~Deligne.
\newblock Th\'eorie de {H}odge. {II}.
\newblock {\em Inst. Hautes \'Etudes Sci. Publ. Math.}, (40):5--57, 1971.

\bibitem{DeligneIII:1971}
P.~Deligne.
\newblock Th\'eorie de {H}odge. {II}.
\newblock {\em Inst. Hautes \'Etudes Sci. Publ. Math.}, (40):5--57, 1971.

\bibitem{Diaconescu:2017}
D.-E. Diaconescu.
\newblock Local curves, wild character varieties, and degenerations.
\newblock {\em Preprint arXiv:1705.05707}, 2017.

\bibitem{Zein-Trang:2014}
F.~El~Zein and L.~D. ung Tr\'ang.
\newblock Mixed {H}odge structures.
\newblock In {\em Hodge theory}, volume~49 of {\em Math. Notes}, pages
  123--216. Princeton Univ. Press, Princeton, NJ, 2014.

\bibitem{Fausk-May:2003}
H.~Fausk, P.~Hu, and J.~P. May.
\newblock Isomorphisms between left and right adjoints.
\newblock {\em Theory Appl. Categ.}, 11:No. 4, 107--131, 2003.

\bibitem{Freed:1994}
D.~S. Freed.
\newblock Higher algebraic structures and quantization.
\newblock {\em Comm. Math. Phys.}, 159(2):343--398, 1994.

\bibitem{Freed-Lurie:2010}
D.~S. Freed, M.~J. Hopkins, J.~Lurie, and C.~Teleman.
\newblock Topological quantum field theories from compact {L}ie groups.
\newblock In {\em A celebration of the mathematical legacy of {R}aoul {B}ott},
  volume~50 of {\em CRM Proc. Lecture Notes}, pages 367--403. Amer. Math. Soc.,
  Providence, RI, 2010.

\bibitem{Freed-Hopkins-Teleman:2010}
D.~S. Freed, M.~J. Hopkins, and C.~Teleman.
\newblock Consistent orientation of moduli spaces.
\newblock In {\em The many facets of geometry}, pages 395--419. Oxford Univ.
  Press, Oxford, 2010.

\bibitem{GP-Gothen-Munoz}
O.~Garc\'ia-Prada, P.~B. Gothen, and V.~Mu\~noz.
\newblock Betti numbers of the moduli space of rank 3 parabolic {H}iggs
  bundles.
\newblock {\em Mem. Amer. Math. Soc.}, 187(879):viii+80, 2007.

\bibitem{GP-Heinloth-Schmitt}
O.~Garc\'ia-Prada, J.~Heinloth, and A.~Schmitt.
\newblock On the motives of moduli of chains and {H}iggs bundles.
\newblock {\em J. Eur. Math. Soc. (JEMS)}, 16(12):2617--2668, 2014.

\bibitem{GP-2019}
{\'A}.~Gonz\'alez-Prieto.
\newblock {$E$}-polynomials of {${\rm SL}_2(\mathbb{C})$}-character varieties
  with semisimple punctures.
\newblock {\em In preparation}.

\bibitem{GP-2018a}
{\'A}.~Gonz\'alez-Prieto.
\newblock Hodge theory of representation varieties via {T}opological {Q}uantum
  {F}ield {T}heories.
\newblock {\em Preprint arXiv:1810.09714}, 2018.

\bibitem{GP-2018b}
{\'A}.~Gonz\'alez-Prieto.
\newblock Stratification of algebraic quotients and character varieties.
\newblock {\em Preprint arXiv:1807.08540}, 2018.

\bibitem{Gonzalez-Prieto:Thesis}
{\'A}.~Gonz\'alez-Prieto.
\newblock Topological {Q}uantum {F}ield {T}heories for character varieties.
\newblock {\em PhD Thesis. Universidad Complutense de Madrid}, 2018.

\bibitem{Gothen}
P.~B. Gothen.
\newblock The {B}etti numbers of the moduli space of stable rank {$3$} {H}iggs
  bundles on a {R}iemann surface.
\newblock {\em Internat. J. Math.}, 5(6):861--875, 1994.

\bibitem{Haugseng}
R.~Haugseng.
\newblock Iterated spans and "classical" topological field theories.
\newblock {\em Preprint arXiv:1409.0837}, 2017.

\bibitem{Hausel:2005}
T.~Hausel.
\newblock Mirror symmetry and {L}anglands duality in the non-abelian {H}odge
  theory of a curve.
\newblock In {\em Geometric methods in algebra and number theory}, volume 235
  of {\em Progr. Math.}, pages 193--217. Birkh\"auser Boston, Boston, MA, 2005.

\bibitem{Hausel-Letellier-Villegas}
T.~Hausel, E.~Letellier, and F.~Rodriguez-Villegas.
\newblock Arithmetic harmonic analysis on character and quiver varieties.
\newblock {\em Duke Math. J.}, 160(2):323--400, 2011.

\bibitem{Hausel-Letellier-Villegas:2013}
T.~Hausel, E.~Letellier, and F.~Rodriguez-Villegas.
\newblock Arithmetic harmonic analysis on character and quiver varieties {II}.
\newblock {\em Adv. Math.}, 234:85--128, 2013.

\bibitem{Hausel-Rodriguez-Villegas:2008}
T.~Hausel and F.~Rodriguez-Villegas.
\newblock Mixed {H}odge polynomials of character varieties.
\newblock {\em Invent. Math.}, 174(3):555--624, 2008.
\newblock With an appendix by Nicholas M. Katz.

\bibitem{Higgins}
P.~J. Higgins.
\newblock Categories and groupoids.
\newblock {\em Repr. Theory Appl. Categ.}, (7):1--178, 2005.
\newblock Reprint of the 1971 original [{{\i}t Notes on categories and
  groupoids}, Van Nostrand Reinhold, London; MR0327946] with a new preface by
  the author.

\bibitem{Hitchin}
N.~J. Hitchin.
\newblock The self-duality equations on a {R}iemann surface.
\newblock {\em Proc. London Math. Soc. (3)}, 55(1):59--126, 1987.

\bibitem{Jardim}
M.~Jardim.
\newblock Nahm transform and spectral curves for doubly-periodic instantons.
\newblock {\em Comm. Math. Phys.}, 225(3):639--668, 2002.

\bibitem{Kassabov-Patotski}
M.~Kassabov and S.~Patotski.
\newblock Character varieties as a tensor product.
\newblock {\em J. Algebra}, 500:569--588, 2018.

\bibitem{Kock:2004}
J.~Kock.
\newblock {\em Frobenius algebras and 2{D} topological quantum field theories},
  volume~59 of {\em London Mathematical Society Student Texts}.
\newblock Cambridge University Press, Cambridge, 2004.

\bibitem{Lauda-Pfeiffer}
A.~D. Lauda and H.~Pfeiffer.
\newblock Open-closed strings: two-dimensional extended {TQFT}s and {F}robenius
  algebras.
\newblock {\em Topology Appl.}, 155(7):623--666, 2008.

\bibitem{Leinster}
T.~Leinster.
\newblock {\em Higher operads, higher categories}, volume 298 of {\em London
  Mathematical Society Lecture Note Series}.
\newblock Cambridge University Press, Cambridge, 2004.

\bibitem{LM}
M.~Logares and V.~Mu\~noz.
\newblock Hodge polynomials of the {$\rm{SL}(2,\Bbb C)$}-character variety of
  an elliptic curve with two marked points.
\newblock {\em Internat. J. Math.}, 25(14):1450125, 22, 2014.

\bibitem{LMN}
M.~Logares, V.~Mu\~noz, and P.~E. Newstead.
\newblock Hodge polynomials of {${\rm SL}(2,\mathbb{C})$}-character varieties
  for curves of small genus.
\newblock {\em Rev. Mat. Complut.}, 26(2):635--703, 2013.

\bibitem{MacLane}
S.~Mac~Lane.
\newblock {\em Categories for the working mathematician}, volume~5 of {\em
  Graduate Texts in Mathematics}.
\newblock Springer-Verlag, New York, second edition, 1998.

\bibitem{Martinez:2017}
J.~Mart\'inez.
\newblock E-polynomials of ${PGL}(2,\mathbb{C})$-character varieties of surface
  groups.
\newblock {\em Preprint arXiv:1705.04649}, 2017.

\bibitem{MM:2016}
J.~Mart\'inez and V.~Mu\~noz.
\newblock E-polynomials of {${SL}(2,\Bbb{C})$}-character varieties of complex
  curves of genus 3.
\newblock {\em Osaka J. Math.}, 53(3):645--681, 2016.

\bibitem{MM}
J.~Mart\'inez and V.~Mu\~noz.
\newblock E-polynomials of the {${\rm SL}(2,\Bbb C)$}-character varieties of
  surface groups.
\newblock {\em Int. Math. Res. Not. IMRN}, (3):926--961, 2016.

\bibitem{Mellit}
A.~Mellit.
\newblock Poincar\'e polynomials of moduli spaces of higgs bundles and
  character varieties (no punctures).
\newblock {\em Preprint arXiv:1707.04214}, 2017.

\bibitem{Mereb}
M.~Mereb.
\newblock On the {$E$}-polynomials of a family of {$SL_n$}-character varieties.
\newblock {\em Math. Ann.}, 363(3-4):857--892, 2015.

\bibitem{Mozgovoy:2012}
S.~Mozgovoy.
\newblock Solutions of the motivic {ADHM} recursion formula.
\newblock {\em Int. Math. Res. Not. IMRN}, (18):4218--4244, 2012.

\bibitem{Nakamoto}
K.~Nakamoto.
\newblock Representation varieties and character varieties.
\newblock {\em Publ. Res. Inst. Math. Sci.}, 36(2):159--189, 2000.

\bibitem{Newstead:1978}
P.~E. Newstead.
\newblock {\em Introduction to moduli problems and orbit spaces}, volume~51 of
  {\em Tata Institute of Fundamental Research Lectures on Mathematics and
  Physics}.
\newblock Tata Institute of Fundamental Research, Bombay; by the Narosa
  Publishing House, New Delhi, 1978.

\bibitem{Peters-Steenbrink:2008}
C.~A.~M. Peters and J.~H.~M. Steenbrink.
\newblock {\em Mixed {H}odge structures}, volume~52 of {\em Ergebnisse der
  Mathematik und ihrer Grenzgebiete. 3. Folge. A Series of Modern Surveys in
  Mathematics [Results in Mathematics and Related Areas. 3rd Series. A Series
  of Modern Surveys in Mathematics]}.
\newblock Springer-Verlag, Berlin, 2008.

\bibitem{Saito:1986}
M.~Saito.
\newblock Mixed {H}odge modules.
\newblock {\em Proc. Japan Acad. Ser. A Math. Sci.}, 62(9):360--363, 1986.

\bibitem{Saito:1989}
M.~Saito.
\newblock Introduction to mixed {H}odge modules.
\newblock {\em Ast\'erisque}, (179-180):10, 145--162, 1989.
\newblock Actes du Colloque de Th\'eorie de Hodge (Luminy, 1987).

\bibitem{Saito:1990}
M.~Saito.
\newblock Mixed {H}odge modules.
\newblock {\em Publ. Res. Inst. Math. Sci.}, 26(2):221--333, 1990.

\bibitem{Schiffmann:2016}
O.~Schiffmann.
\newblock Indecomposable vector bundles and stable {H}iggs bundles over smooth
  projective curves.
\newblock {\em Ann. of Math. (2)}, 183(1):297--362, 2016.

\bibitem{Schurmann:2011}
J.~Sch\"urmann.
\newblock Characteristic classes of mixed {H}odge modules.
\newblock In {\em Topology of stratified spaces}, volume~58 of {\em Math. Sci.
  Res. Inst. Publ.}, pages 419--470. Cambridge Univ. Press, Cambridge, 2011.

\bibitem{Mozgovoy-Schiffmann}
M.~Sergey and O.~Schiffmann.
\newblock Counting higgs bundles and type a quiver bundles.
\newblock {\em Preprint arXiv:1705.04849}, 2017.

\bibitem{Schulman}
M.~Shulman.
\newblock Framed bicategories and monoidal fibrations.
\newblock {\em Theory Appl. Categ.}, 20:No. 18, 650--738, 2008.

\bibitem{Simpson:parabolic}
C.~T. Simpson.
\newblock Harmonic bundles on noncompact curves.
\newblock {\em J. Amer. Math. Soc.}, 3(3):713--770, 1990.

\bibitem{Simpson:1992}
C.~T. Simpson.
\newblock Higgs bundles and local systems.
\newblock {\em Inst. Hautes \'Etudes Sci. Publ. Math.}, (75):5--95, 1992.

\bibitem{SimpsonI}
C.~T. Simpson.
\newblock Moduli of representations of the fundamental group of a smooth
  projective variety. {I}.
\newblock {\em Inst. Hautes \'Etudes Sci. Publ. Math.}, (79):47--129, 1994.

\bibitem{SimpsonII}
C.~T. Simpson.
\newblock Moduli of representations of the fundamental group of a smooth
  projective variety. {II}.
\newblock {\em Inst. Hautes \'Etudes Sci. Publ. Math.}, (80):5--79 (1995),
  1994.

\bibitem{Strominger-Yau-Zaslow}
A.~Strominger, S.-T. Yau, and E.~Zaslow.
\newblock Mirror symmetry is {$T$}-duality.
\newblock {\em Nuclear Phys. B}, 479(1-2):243--259, 1996.

\end{thebibliography}
\bibliographystyle{abbrv}

\end{document}